[1994/12/01]
\documentclass[reqno,13pt]{amsart}
\usepackage{bbm}
\def\1{\mathbbm{1}}

\usepackage{amsmath,amsthm}
\usepackage[mathscr]{euscript}
\usepackage{graphicx}
\usepackage{amsmath,amssymb}
\usepackage{circuitikz}

\usepackage{pstricks}
\usepackage{pst-node}

\usepackage{tikz}
\usetikzlibrary{arrows,chains,matrix,positioning,scopes}
\tikzstyle{element}=[rectangle,draw,fill=white, line width=1pt]
\tikzstyle{terminal}=[circle,draw, scale=0.3, line width=1pt,red]
\tikzstyle{fleche}=[->,>=stealth', very thick]
\tikzstyle{fleche1}=[->,>=stealth', very thick, red]
\usepackage{placeins}
\definecolor{darkgoldenrod4}{rgb}{0.55,0.4,0.55}
\definecolor{maroon4}{rgb}{0.55,0.11,0.38}
\definecolor{indianred}{rgb}{0.8,0.36,0.36}
\definecolor{purple1}{rgb}{0.61,0.19,1}
\definecolor{goldenrod1}{rgb}{1,0.76,0.15}
\definecolor{indianred3}{rgb}{0.8,0.33,0.33}
\definecolor{red4}{rgb}{0.55,0,0}
\definecolor{darkslategray}{rgb}{0.18,0.31,0.31}
\definecolor{firebrick}{rgb}{0.7,0.13,0.13}
\definecolor{slateblue3}{rgb}{0.41,0.35,0.8}
\definecolor{mediumorchid4}{rgb}{0.48,0.22,0.55}
\definecolor{thistle4}{rgb}{0.55,0.48,0.55}
\definecolor{rltred}{rgb}{0.75,0,0}
\definecolor{rltgreen}{rgb}{0,0.5,0}
\definecolor{oneblue}{rgb}{0,0,0.75}
\definecolor{marron}{rgb}{0.64,0.16,0.16}
\definecolor{forestgreen}{rgb}{0.13,0.54,0.13}
\definecolor{purple}{rgb}{0.62,0.12,0.94}
\definecolor{dockerblue}{rgb}{0.11,0.56,0.98}
\definecolor{freeblue}{rgb}{0.25,0.41,0.88}
\definecolor{myblue}{rgb}{0,0.2,0.4}

\def\R{{\mathbb{R}}}

\def\al{\alpha}

\def\ga{\gamma}

\def\la{\lambda}

\def\t{\tau}

\def\calM{{\mathcal{M}}}

\def\calL{{\mathcal{L}}}

\def\calS{{\mathcal{S}}}

\def\R{\mathbb R}
\def\calN{\mathcal{N}}

\def\C{\mathbb C}

\def\N{\mathbb N}

\def\F{\mathbb F}

\def\T{\mathbb T}

\newtheorem{definition}{Definition}[section]

\newtheorem{lemma}{Lemma}[section]
\newtheorem{proposition}{Proposition}[section]
\newtheorem{corollary}{Corollary}[section]
\newtheorem{theorem}{Theorem}[section]
\theoremstyle{remark}
\newtheorem{remark}{Remark}[section]

\title{Positivity of infinite-dimensional linear systems} 
\thanks{Department of Mathematics, Zhejiang Normal University, Jinhua, Zhejiang 321004, P.R. China}
\author{Yassine El Gantouh} 


\begin{document}
	\maketitle
	
	\renewcommand{\sectionmark}[1]{}
	\begin{abstract}
			In this paper, we investigate the well-posedness and positivity property of infinite-dimensional linear system with unbounded input and output operators. In particular, we characterize the internal and external positivity for this class of systems. This latter effort is motivated in part by a complete description of well-posed positive control/observation systems. An interesting feature of positive well-posed linear systems is that weak regularity and strong regularity, in the sense of Salamon and Weiss, are equivalent. Moreover, we provide sufficient conditions for zero-class admissibility for positive semigroups. In the context of positive perturbations of positive semigroups, we establish two perturbation results, namely the Desch-Schappacher perturbation and the Staffans-Weiss perturbation. As for illustration, these findings are applied to investigate the existence and uniqueness of a positive mild solution of the linear Boltzmann equation with non-local boundary conditions on finite network.
	\end{abstract}

	\section{Introduction}
\label{intro}
We aim at studying the well-posedness and positivity properties of solutions of linear time-invariant systems on infinite dimensional spaces, which are described by the state-space representation 
\begin{subequations}
	\begin{align}
		\dot{z}(t)&= Az(t)+Bu(t), \qquad
		t > 0,\label{eq:1} \\ 
		z(0)&=x
		\nonumber, \\
		y(t)&= Cz(t)+Du(t), \qquad t > 0,\label{eq:2}
	\end{align}
\end{subequations}
where $ z$ is the state, $x$ is a positive initial data, $u$ is a positive input, $y$ is the output, $ A $ generates a strongly continuous semigroup, $D$ is a bounded linear operator, and $ B,C $ are unbounded linear operators on appropriate Banach spaces. The unboundedness of the operators $B$ and $C$ leads to the notion of admissibility in order to have continuous dependence of $L^p$-controls or continuous dependence of $L^p$-observations on states \cite{Salam,TW1,WO,WC}. If one has in addition continuous dependence of $L^p$-observations on $L^p$-controls, then the input-output system \eqref{eq:1}-\eqref{eq:2} is called an $L^p$-well-posed linear system \cite{Salam,Staf}. This class of systems, known as \emph{Salamon-Weiss systems}, covers, in particular, PDEs with boundary control and point observation systems, as well as FDEs with delays in the input and output variables, see e.g., \cite{Salam,Staf,TW1}.

Many mathematical models can be written in the abstract formulation \eqref{eq:1}-\eqref{eq:2}. This fact has motivated several works to investigate the qualitative and quantitative properties of $L^p$-well-posed linear systems (see e.g., \cite{CM,HMR,Salam,Sc1,WR,WCur,Zwart}). For instance, in \cite{WR} Weiss introduced an important subclass of $L^p$-well-posed linear system, called \emph{regular systems}, where one can also find an elegant feedback theory. In \cite{WCur} the authors studied the concept of stabilizability of $L^p$-well-posed and regular linear systems. The linear quadratic control problem for $L^p$-well-posed linear systems has been addressed in several works in the literature (see e.g., \cite{Salam,Staf} and references therein). In the context of perturbations theory of one-parameter operator semigroups, it has been shown that the feedback theory of $L^p$-well-posed and regular linear systems provide a useful analytical approach when dealing with boundary perturbations \cite{CM,HMR,Salam,Zwart}. A generalization of $L^p$-well-posed and regular linear system to the non-autonomous setting is given in Schnaubelt \cite{Sc1}. Further discussion on the practical interest and nearly complete overview of this field can be found in \cite{TW2}.

In many practical applications some constraints on the state and/or on the control need to be imposed.  This is for instance the case of evolutionary networks (heat conduction, transportation networks, etc.) where realistic models have to take into account that the state represents some physical quantity which must necessarily remain positive. Theses type of systems are generally referred to as \emph{positive systems} \cite{FR,Luenberger}. As for ordinary differential equations, in the finite dimensional situation one has a well established theory of positive systems, see e.g., \cite{FR} or \cite{RV} for a comprehensive survey and detailed lists of references. However, its infinite-dimensional counterpart is much less present in the literature, due to the fact that finite-dimensional results cannot be straightforwardly generalized (see e.g., \cite{AADHW,Tilman}). For instance, when working in such a framework, it is not clear a priori if the system \eqref{eq:1}-\eqref{eq:2} has a positive solution.

In this paper, we study the input-output system \eqref{eq:1}-\eqref{eq:2}, in the Banach lattice setting, with a particular aim to explore new techniques and new questions for infinite-dimensional positive linear systems. To this end, we present a new and very operator theoretical approach which allows us to
achieve the following three goals:
\begin{itemize}
	\item[(1)]  We describe the structural properties of $L^p$-admissible positive control and observation operators.
	\item[(2)]  We characterize positive $L^p$-well-posed regular linear systems.
	\item[(3)]  We give new insight into the perturbation theory, namely, positive perturbations of positive semigroups.
\end{itemize}
Let us elaborate on these points in more detail.

(1) As a first contribution, we describe the structural properties of $L^p$-admissible positive control/observation operators in the Banach lattice setting. In fact, we show that the admissibility of positive observation and control operators is fully determined by their behavior on the positive part of $D(A)$ and with respect to the positive inputs (Lemma \ref{S2.L6} and Remark \ref{S2.R2}), respectively. Then, by a slight modification of the proof of \cite[Lemma 2.1 ]{Voigt}, in Lemma \ref{S2.P3} we show that positive observation operators are $L^1$-admissible provided that the output-space is an $AL$-space. Further, we provide sufficient conditions with respect to the resolvent operator for zero-class admissibility of positive control and observation operators (Theorem \ref{S2.T1} and Proposition \ref{S2.P4}).

(2)  We characterize the internal and external positivity of the input-output system \eqref{eq:1}-\eqref{eq:2}. The latter effort is motivated in part by a complete description of $ L^p $-admissible positive control/observation operators. An interesting feature of positive $L^p$-well-posed linear systems is that weak regularity and regularity are equivalent, see Theorem \ref{S4.P4}. This fact implies that the closed-loop system obtained from a positive regular linear system with a positive admissible feedback operator exhibits all properties of the input-output system \eqref{eq:1}-\eqref{eq:2} (Theorem \ref{S4.T3} and Remark \ref{Remark}).

(3) As a standard application of the developed theory, we establish two results on positive perturbations of positive semigroups (Theorem \ref{S2.P2} and Theorem \ref{S4.CB}), namely the Desch-Schappacher perturbation and the Staffans-Weiss perturbation. The proof of Theorem \ref{S2.P2} is based on a generation result of Arendt \cite[Theorem 2.5]{Arendt}, the origin of which can be traced back to Batty and Robinson \cite{BaRo}. Theorem \ref{S4.CB} is obtained as a consequence of Theorem \ref{S4.T3}.

We now depict some related works in the literature. To the best of our knowledge, positive $L^p$-well-posed and regular linear system were first studied in \cite{Wint}. In this thesis, the author introduced the class of positive $L^p$-well-posed linear systems using an abstract framework. There, the author also observed that on finite-dimensional control spaces, positive control operators are zero-class admissible for positive semigroups. There are not other works in the literature studying positive $L^p$-well-posed linear systems. It seems that the attention was mainly directed towards positive perturbations of positive semigroups. For instance, in \cite{BJVW} the authors showed a generation result for perturbations of Desch-Schappacher type. In \cite{BBH}, the authors addressed the positivity of perturbed boundary value problems and also they showed a stability result of a class of linear hyperbolic systems with delay at the boundary. An important recent contribution is presented in \cite{GM}, in which the authors developed uniform small-gain conditions to efficiently characterize the exponential stability of positive linear  discrete-time systems in ordered Banach spaces. More recently, using the concept of positive $L^p$-well-posed and regular linear systems introduced in this paper, in \cite{El1} we provided necessary and sufficient conditions for controllability under positivity constraints of non-homogeneous boundary value control problems.

The paper is organized as follows. In the remaining part of this introduction we recall some well-known results on Banach lattices theory and positive semigroups. In Section \ref{Sec2} we investigate the positivity properties of solutions of linear control systems. In Section \ref{Sec3}, we study the admissibility properties of positive observation operators. Section \ref{Sec4} deals with the positivity property of the input-output system \eqref{eq:1}-\eqref{eq:2} and the corresponding closed-loop systems. Finally, in Section \ref{Sec5} we establish the existence and uniqueness of a positive mild solution of the transport equation on finite network with scattering at vertices.

\textbf{Notations and preliminaries.} Before we come to the description of our main results, we need to fix some notation and recall some definition on \emph{vector lattices}- also called \emph{Riesz spaces}. The reader is referred to the books \cite{CHARALAMBOS} by C.D. Aliprantis and O. Burkinshaw or \cite{Schaf} by H.H. Schaefer, which were our main references, if more details are required. Let $E$ be a real vector space and $\leq$ be a partial order on this space. Then $E$ is said to be an ordered vector space if it satisfies the following properties:
\begin{itemize}
	\item[$(\emph{i})$] if $f,g\in E$ and $f\leq g$, then  $f+h\leq g+h$ for all $h\in E$.
	\item[$(\emph{ii})$] if $f,g\in E$ and $f\leq g$, then  $\alpha f\leq \alpha g$ for all $\alpha\geq 0$.
\end{itemize}
If, in addition, $E$ is lattice with respect to the partial ordered, that is, $\sup\{f,g\}$ and $\inf\{f,g\}$ exist for all $f,g\in E$, then $E$ is said to be a \emph{vector lattice}. For an element $f$ of a \emph{vector lattice} $E$, the \emph{positive part} of $f$ is defined by $f_+:=\sup\{f,0\}$, the \emph{negative part} of $f$ by $f_-:=\sup\{-f,0\}$ and the absolute value of $f$ by $\vert f\vert:=\sup\{f,-f\}$, where $0$ is the zero element of $E$. The set of all positive elements of $E$, denoted by $E_+$, is a convex cone with vertex $0$. In particular, it generates a canonical ordering $ \leq $ on $ E $ which is given by: $f\leq g$ if and only if $g-f\in E_+$. Note that the positive cone of a normed vector lattice is closed. 
A norm complete vector lattice $ E $ such that its norm satisfies the following property
\begin{displaymath}
	\vert f\vert \leq \vert g\vert \quad \implies\quad \Vert f\Vert \leq \Vert g\Vert
\end{displaymath}
for $f,g\in X$, is called \emph{Banach lattice}. If $E$ is a Banach lattice, its topological dual $E'$, endowed with the dual norm and the dual order, is also a Banach lattice. A norm on a Banach lattice is \emph{order continuous} if for each generalized sequence $(x_\alpha)$ such that $x_\alpha \downarrow 0$ in $E$, the sequence $(x_\alpha)$ converges to $0$ for the norm $\Vert\cdot\Vert$, where the notation $x_\alpha \downarrow 0$ means that the sequence $(x_\alpha)$ is decreasing (in symbol $\downarrow$), its infimum exists and $\inf (x_\alpha)=0$ (see, e.g., \cite[Definition 4.7]{CHARALAMBOS}). 

Now, if we denote by $ \mathcal{L}(E,F) $ the Banach algebra of all linear bounded operators from an order Banach space $E$ to an order Banach space $F$, then an operator $ P\in \mathcal{L}(E,F)  $ is positive if $ P E_+\subset F_+ $. An everywhere defined positive operator from a Banach lattice to a normed vector lattice is bounded, see e.g., \cite[Theorem II.5.3]{Schaf}. Note that for $P\in \mathcal{L}_+(E,F)$ we have
\begin{align*}
	\Vert P\Vert= \sup_{x\ge 0,\Vert x\Vert\le 1} \Vert Px\Vert,
\end{align*}
see e.g., \cite[Proposition 10.22]{BFR}-a). The set of all positive operators from a Banach lattice $E$ to another Banach lattice $F$, denoted by $ \mathcal{L}_+(E,F) $, is a convex cone in the space $ \mathcal{L}(E,F) $. 

Let $X$ be a Banach lattice and $ (A,D(A)) $ be the generator of a C$_{0} $-semigroup $ \mathbb{T}:=(T(t))_{t\geq 0} $ on $ X $. The type of $ \mathbb{T} $ is defined by $ \omega_{0}(A):=\inf \lbrace t^{-1}\log \Vert T(t)\Vert: t> 0\rbrace $. We denote by $ \rho (A) $ the resolvent set of $ A $, i.e., the set of all $ \mu \in \C $ such that $ \mu I_X -A $ has an inverse in $\calL(X)$ with $I_X$ denote the identity operator in $X$. By $ R(\mu,A): =(\mu I_X-A)^{-1} $ we denote the resolvent operator of $ A $. The complement of $ \rho (A) $, is called the spectrum and is denoted by $\sigma(A):= \C\backslash \rho (A)$. The so-called spectral radius of $A$ is defined by $r(A):=\sup\{\vert \mu\vert :\; \mu \in \sigma(A)\}$. Also, recall that the spectral bound $s(A)$ of $A$ is defined by $s(A):=\sup\{{\rm Re}\, \mu :\; \mu \in \sigma(A)\}$. A linear operator $A$ on a Banach lattice $X$ is called resolvent positive if there exists $\omega\in \mathbb{R}$ such that $(\omega,\infty)\subset \rho (A) $ and $R(\mu,A)\geq 0$ for each $\mu> \omega$. It follows from \cite[Corollary 2.3]{Arendt} that a C$_{0} $-semigroups on a Banach lattice is positive if and only if the corresponding generator $A$ is resolvent positive. On the other hand, by $ X_{1} $ we denote the order Banach space $ D(A) $ endowed with the norm $ \Vert x\Vert_{1}:=\Vert (\mu I_X-A)x\Vert $ for some $ \mu \in \rho (A) $. The extrapolation space associated with $ X $ and $ A $, denoted by $ X_{-1} $, is the completion of $ X $ with respect to the norm  $\Vert x\Vert_{-1}:=\Vert R(\mu , A)x\Vert $ for $ x\in X $ and some $ \mu \in \rho (A) $. Note that the choice of $ \mu  $ is not important, since by the resolvent equation different choices lead to equivalent norms on $ X_{1} $ and $ X_{-1} $. Note also that $f\in X_{-1}$ is positive if $f$ belongs to the closure of $X_+$ in $X_{-1}$. In addition, we have    
$$
X_{1}\subset X\subset X_{-1}\quad {\rm and}\quad X_+=X\cap X_{-1,+},$$ where $X_{-1,+}$ denotes the closure of $X_+$ in $X_{-1}$ (see \cite[Proposition 2.3]{BJVW} for the second statement). We point out that the extrapolation space $X_{-1}$ is not, in general, a Banach lattice (more precisely, a vector lattice), cf. \cite[Remark 2.5]{BJVW}. The unique extension of $ \mathbb{T} $ on $ X_{-1} $ is a $C_{0}$-semigroup which we denote by $ \mathbb{T}_{-1} $ and whose generator is denoted by $ A_{-1} $. For an overview of the theory of positive C$_0$-semigroups, we refer for instance to the monographs \cite{BFR,Nagel}. 

Let $ U$ be a Banach lattice, $\R_+:=[0,+\infty)$, $p\in [0,+\infty)$ and $ \alpha>s(A)$. We denote by $L^{p}_{\alpha}(\R_+;U)$ the space of all functions of the form $ v(t)=e^{\alpha t}u(t) $, where $ u\in L^{p}(\R_+;U) $. Moreover, let $L_{loc,+}^p(\R_+;U)$ denote the set of positive control functions $u$ in $L_{loc}^p(\R_+;U)$ such that $u(t)\in U_+$ almost everywhere in $\R_+$, where we regard $L_{loc}^p(\R_+;U)$ as a lattice ordered Fr\'{e}chet space with the seminorms being the $L^p$ norms on the intervals $[0,n]$, $n\in\N$. For $\t\ge 0$, the truncation operator $\mathcal{P}_\t$ and the left shift operator $\calS_\t$ are defined by
\begin{align*}
	(\mathcal{P}_\t u)(t)=
	\begin{cases}
		u(t),&  {\rm for }\;\; t\in [0,\t),\cr
		0,& {\rm for }\;\; t\ge \t,
	\end{cases},\qquad 
	(\mathcal{S}_\t u)(t)=
	\begin{cases}
		0,&  {\rm for }\;\; t\in [0,\t),\cr
		u(t+\t),& {\rm for }\;\; t\ge \t,
	\end{cases}
\end{align*}
respectively.

	\section{Admissible positive control operators} \label{Sec2}
Let $X$ and $U$ be Banach lattices. The main concern of this section is the positivity and well-posedness, in the Banach lattice setting, of the abstract differential equation
\begin{equation}\label{S2.1}
	\dot{z}(t) = A_{-1} z(t)+B u(t),\qquad t> 0,\qquad z(0)=x,
\end{equation}
where $ A$ generates a positive C$_0$-semigroup $ \mathbb{T}$ on $ X $ and $B$ is a bounded linear operator from $U$ into $ X_{-1} $. An integral solution of \eqref{S2.1} is given by 
\begin{equation}\label{S2.2}
	z(t;x,u)=T(t)x+\int_{0}^{t}T_{-1}(t-s)Bu(s)ds, \; \; t\geq 0,
\end{equation}
for $ x\in X $ and $u\in  L^{p}_{loc}(\R_+;U) $, where the integral in \eqref{S2.2} is calculated in $ X_{-1} $. However from the perspective of control theory, we usually look for continuous $ X $-valued functions $z(.;x,u)$. This motivates the following definition, see \cite{WC}.
\begin{definition}\label{S2.D1}
	An operator $ B\in \mathcal{L}(U,X_{-1}) $ is called $ L^p $-admissible control operator for $A$ (with $p\in [1,+\infty)$) if, for some $\tau >0 $, 
	\begin{equation}\label{S2.2'}
		\Phi_{\tau}u:=\int_{0}^{\tau}T_{-1}(\tau-s)Bu(s)ds,
	\end{equation}
	takes values in $X$ for any $ u\in L^{p}(\R_+;U) $.
\end{definition}
Note that, if $ B $ is $ L^p $-admissible for some $ \t>0 $, then
\begin{equation*}
	\Vert \Phi_{\tau}u\Vert \leq \kappa\Vert u\Vert_{ L^{p}([0,\tau];U)},
\end{equation*}
for all $ u\in L^{p}([0,\tau];U) $ and a constant $\kappa:=\kappa(\t)>0 $. Note also that, if $ B $ is $ L^p $-admissible for some $ \t>0 $, then for every $\tau\ge 0$ we have 
\begin{align*}
	\Phi_\tau\in \calL(L^{p}(\R_+;U),X).
\end{align*}
In addition, if $ \underset{\t\mapsto 0 }{\lim} \kappa(\t)=0$, then we say that $B$ is a zero-class $L^p$-admissible control operator for $A$. We denote by $ \mathfrak{B}_p(U,X,\mathbb{T}) $ the vector space of all $L^p$-admissible control operators $B$, which is a Banach space with the norm 
\begin{displaymath}
	\Vert B\Vert_{\mathfrak{B}_p}:=\sup_{\Vert u\Vert_{_{L^{p}([0,\t];U)}}\leq 1}\left\Vert \int_{0}^{\tau}T_{-1}(\tau-s)Bu(s)ds\right\Vert,
\end{displaymath} 
where $\t> 0$ is fixed, see \cite{WC} for more details.

We have the following characterization of the well-posedness of \eqref{S2.1}.
\begin{proposition}\label{S2.P1}
	Let $ X , U $ be Banach lattices, let $ A $ generates a strongly continuous semigroup $ \mathbb{T}$ on $ X $ and $ B\in \mathcal{L}(U,X_{-1}) $. Then, the following assertions are equivalent:
	\begin{itemize}
		\item[(\emph{i})] For any $x\in X_+$ and $ u\in L^{p}_+(\R_+;U) $, the solution $ z(\cdot) $ of \eqref{S2.1} remains in $ X_+ $.
		\item[(\emph{ii})] $\mathbb{T},B $ are positive and for some $ \tau>0,\,  \Phi_\tau u \in X_+ $ for all $ u\in L^{p}_+([0,\tau];U) $. 
	\end{itemize}
	In addition, if one of these conditions is satisfied then the differential equation \eqref{S2.1} has a unique positive mild solution $ z(\cdot)\in C(\R_+;X)$ given by $z(t)=T(t)x+\Phi_tu$ for all $t\ge 0$ and $(x,u)\in X\times L^p_{loc}(\R_+,U)$.
\end{proposition} 	
\begin{proof}
	The equivalence $ (\emph{i}) \Leftrightarrow (\emph{ii}) $ is straightforward from \eqref{S2.2}. By $ (\emph{ii}) $ we have $ \Phi_\t u \in X_+ $ for some $\t>0$ and for every $u\in L^{p}_+([0,\tau];U)$. Thus, by virtue of \cite[Theorem 1.10]{CHARALAMBOS}, $ \Phi_\t  $ has a unique extension to a positive operator from $ L^{p}([0,\tau];U) $ to $ X $, since $ \Phi_\t $ is linear. Thus, according to Definition \ref{S2.D1}, $ B $ is an $L^p$-admissible control operator for $ A $ and hence $\Phi_\tau\in \calL(L^{p}(\R_+;U),X)$ for any $\tau\ge 0$. Therefore, the abstract differential equation \eqref{S2.1} has a unique mild solution, cf. \cite[Theorem 4.3.1]{Staf}.
\end{proof}

\begin{remark}\label{S2.R1}
	Using linear control systems terminology (see e.g., \cite[Definition 2.1]{WC}) and according to \cite[Theorem 1.10]{CHARALAMBOS}, one can see that abstract positive control systems are fully determined by their behaviour on the positive cone. In fact, let $ \mathbb{T} $ be a positive $ C_0 $-semigroup on $ X $, and let $ \Phi:=(\Phi_\t)_{\t\geq 0} $ be an additive family of operators $ \Phi_\t: L^{p}_+(\R_+;U)\to X_+ $ satisfying  
	\begin{align}\label{S2.4}
		\Phi_{\t+t}u=T(\t)\Phi_t \mathcal{P}_t u+\Phi_\t\mathcal{S}_t u,
	\end{align}
	for $ \t,t\geq 0 $ and $ u\in L^{p}_+([0,\t+t],U)$. Thus, according to \cite[Theorem 1.10]{CHARALAMBOS}, the map $ \Phi_\t $ has a unique extension to a positive linear operator from $ L^{p}(\R_+;U) $ into $ X $ for all $\t\ge 0$. Moreover, the extension, denoted by $\Phi_\t $ again, is given by
	\begin{align}\label{composition}
		\Phi_{\t}u= \Phi_{\t}u_+-\Phi_{\t}u_-,
	\end{align}
	for all $ u\in L^{p}(\R_+,U)$ and $\t\ge 0$, cf. \cite[Theorem 1.10]{CHARALAMBOS}. Thus, for every $\t\ge 0$, $ \Phi_\t\in \mathcal{L}(L^{p}(\R_+;U),X) $ as an everywhere defined positive operator is bounded, cf. \cite[Theorem II.5.3]{Schaf}. The equation \eqref{composition} further yields that $\Phi_\t$ satisfies \eqref{S2.4} for arbitrary $ u\in L^{p}(\R_+;U) $. Therefore, $(\mathbb{T},\Phi)$ is a positive linear control system on $X$ and $L^{p}(\R_+;U)$. In addition, according to \cite[Theorem 3.9]{WC}, there exists a unique positive control operator $ B\in \mathcal{L}(U,X_{-1}) $ such that $ \Phi_\t $ is given by the integral representation formula \eqref{S2.2'}. Note that the positivity of $ B $ follows from the closedness of $ X_{-1,+} $, since $ Bv=\underset{\t\mapsto 0 }{\lim}  \tfrac{1}{\t}\Phi_{\t}v $ (in $ X_{-1} $) for all $ v\in U $.
\end{remark}
Next, according to \cite[Theorem 5.3.1]{ACMF}, $ u $ and $ z $ from \eqref{S2.2} have Laplace transforms related by	
\begin{align}\label{laplace}
	\hat{z}(\mu)=R(\mu,A)x+\widehat{(\Phi_\cdot u)}(\mu),\quad \text{ with }\quad \widehat{(\Phi_\cdot u)}(\mu)=R(\mu,A_{-1})B \hat{u}(\mu) ,
\end{align}
for all $\Re e\,\mu> \alpha$, where $\hat{u}$ denote the Laplace transform of $u$.
\begin{lemma}\label{S2.LB}
	Let $ X,U$ be Banach lattices and let $ \mathbb{T} $ be a positive $ C_0 $-semigroup on $ X $. Then, the following assertions are equivalent:
	\begin{itemize}
		\item[$(\emph{i})$] $B$ is positive.
		\item[$(\emph{ii})$] $R(\mu,A_{-1})B$ is positive for all $\mu > s(A)$.
		\item[$(\emph{iii})$] $\Phi_t$ is positive for all $t\geq 0$.
	\end{itemize}
\end{lemma}	
\begin{proof}
	The equivalence of (\emph{i}) and (\emph{ii}) is shown in \cite[Remark 2.2]{BJVW}. On the other hand, the integral formula \eqref{S2.2'} together with the fact that
	$
	B v=\underset{t\mapsto 0 }{\lim}\frac{1}{t}\Phi_t v
	$
	yield the equivalence of $(\emph{i})$ and $ (\emph{iii})$ (as $ X_+= X_{-1,+}\cap X$). 
\end{proof}
\begin{remark}\label{S2.R2}
	Looking at Lemma \ref{S2.LB} and according to \cite[Theorem 1.10]{CHARALAMBOS}, one could relax the definition of $L^p$-admissibility of positive control operators for positive semigroups to, for some $\t>0$, 
	\begin{displaymath}
		\Phi_\t u\in  X_+,\qquad  \forall\, u\in L^p_+(\R_+;U).
	\end{displaymath}
	Thus, $ B $ is an $ L^p $-admissible positive control operator for $A$ if and only if for $\tau>0$ there exists a constant $\kappa:=\kappa(\t)>0 $ such that for every step function $u:[0,\tau]\to U_+$,
	\begin{equation}\label{positive-input}
		\Vert \Phi_{\tau}u\Vert \leq \kappa\Vert u\Vert_{ L^{p}([0,\tau];U)},
	\end{equation}
	In fact, the positivity of the input-map $ \Phi_{\tau} $ implies that $ \vert\Phi_{\tau}u\vert\leq \Phi_{\tau}\vert u\vert  $, for arbitrary $ U $-valued step function $ u $. It follows that
	$$
	\Vert \Phi_{\tau}u\Vert \leq \Vert \Phi_{\tau}\vert u\vert\Vert \leq  \kappa\Vert \vert u\vert\Vert_{ _{L^{p}([0,\t];U)}}=\kappa\Vert u\Vert_{ _{L^{p}([0,\t];U)}},
	$$
	and hence, according to \cite[Remark 6.7]{WO}, $ B $ is an $ L^p $-admissible positive control operator for $A$.
\end{remark}
\begin{remark}
	It is worth to note that for an $L^p$-admissible positive control operator $B$, we have 
	\begin{align*}
		\Vert B\Vert_{\mathfrak{B}_p}=\sup_{u\ge 0,\Vert u\Vert\leq 1}\left\Vert \int_{0}^{\tau}T_{-1}(\tau-s)Bu(s)ds\right\Vert,
	\end{align*}
	for some $\tau>0$. So, if we denote by $ \mathfrak{B}_{p,+}(U,X,\mathbb{T}) $ the set of all $L^p$-admissible positive control operators, then it is not difficult to see that $ \mathfrak{B}_{p,+}(U,X,\mathbb{T}) $ is the positive convex cone of $ \mathfrak{B}_{p}(U,X,\mathbb{T}) $. Therefore, it generates the order: for $ B,B'\in  \mathfrak{B}_p(U,X,\mathbb{T})$, we have $ B'\leq B$ whenever $ B-B'\in \mathfrak{B}_{p,+}(U,X,\mathbb{T}) $.
\end{remark}

The next result provide a sufficient condition ensuring the admissibility of positive control operators for positive semigroups. 
\begin{theorem}\label{S2.T1}
	Let $ X , U$ be Banach lattices and let $ A$ be a densely defined resolvent positive operator such that
	\begin{align}\label{S2.5}
		\Vert R(\mu_0,A)x\Vert\geq c \Vert x \Vert, \qquad x\in X_+,
	\end{align}
	for some $ \mu_0>s(A) $ and constant $c >0$. Let $ B\in \mathcal{L}(U,X_{-1}) $ be a positive control operator. Then $ B $ is an $ L^p $-admissible positive control for $ A $.
\end{theorem}
\begin{proof}
	Note that the inverse estimate \eqref{S2.5} implies that $A$ generates a positive $ C_0 $-semigroup $\mathbb{T}$ on $ X $, cf. \cite[Theorem 2.5]{Arendt}. Since step functions are dense in $ L^{p} $ for $ 1\leq p <\infty $, then it suffices to show that, for every step function $ v:[0,\t]\to  U $, the estimate \eqref{positive-input} holds . In fact, let $ v  $ be a positive step function i.e., $ v(t)=\sum_{k=1}^{n} v_k \1_{[\t_{k-1},\t_{k}[}(t) $ where $ v_1,\ldots,v_n\in U_+ $, and $ \1_{[\t_{k-1},\t_{k}[} $ denotes the indicator function of $ [\t_{k-1},\t_{k}[$ for $ 0=\t_0<\ldots< \t_n=\t $. Clearly, for any positive step function $ v $ we have $ \Phi_{\tau}v \in X_+$. Using the inverse estimate \eqref{S2.5}, we obtain
	\begin{align*}
		c\Vert \Phi_{\tau}v\Vert& \le \left\Vert R(\mu_0,A)  \Phi_{\tau}v\right\Vert\\
		&= \left\Vert  \sum_{k=1}^{n}\int_{\t_{k-1}}^{\t_k}T(\tau -s)R(\mu_0,A_{-1})B v_kds\right\Vert\\
		&\leq \int_{0}^{\t}\Vert T(\tau -s)R(\mu_0,A_{-1})Bv(s)\Vert ds\\
		&\leq  M\Vert R(\mu_0,A_{-1})B\Vert  \t^{\tfrac{1}{q}}e^{w\t} \Vert v\Vert_{ _{L^{p}([0,\t];U)}},
	\end{align*}
	where we have used the estimate $ \Vert T(t)\Vert \leq Me^{wt} $ for $ w>s(A) $ and the H\"{o}lder's inequality (for $ \frac{1}{p}+\frac{1}{q}=1 $). Thus,
	\begin{align*}
		\Vert \Phi_{\tau}v\Vert \le \frac{1}{c}M\Vert R(\mu_0,A_{-1})B\Vert  \t^{\tfrac{1}{q}}e^{w\t} \Vert v\Vert_{ _{L^{p}([0,\t];U)}},
	\end{align*}
	for all positive step function $v$. Hence, according to Remark \ref{S2.R2}, $B$ is an $L^p$-admissible positive control operator for $A$. This completes the proof.
\end{proof}

\begin{remark}\label{S2.R4}
	In the proof of the above proposition we have used the fact that $ R(\mu_0,A_{-1})B $ is positive as an operator in $ \mathcal{L}(U,X) $ for $ \mu_0>s(A) $, see \cite[Proposition 1.7.2]{BaRo}. The proof of Theorem \ref{S2.T1} further yields that $ B $ is a zero-class $L^p$-admissible positive control operator for all $p>1$, since $$ \lim_{\t\mapsto 0 } \kappa(\t)=\lim_{\t\mapsto 0 }\frac{1}{c}M\Vert R(\mu_0,A_{-1})B\Vert  \t^{\tfrac{1}{q}}e^{w\t}=0, \quad q<+\infty.$$ Note that the hypothesis \eqref{S2.5} requires an inverse estimate with respect to the Hille-Yosida Theorem. Note also that, in practical situations, the Banach lattice $X$ is an $L^1$-space (An example satisfying \eqref{S2.5} is given in Section \ref{Sec5}).
\end{remark}

With the help of the above theorem we obtain the following perturbation result of positive semigroups.
\begin{theorem}\label{S2.P2}
	Let $A$ be a densely defined resolvent positive operator and for any $\mu\ge \mu_0$, for some $\mu_0>s(A)$, there exists $c>0$ such that
	\begin{align}\label{inverse-second}
		\Vert R(\mu,A)x\Vert\geq c \Vert x \Vert, \qquad \forall x\in X_+.
	\end{align}
	If $ B\in \mathcal{L}(X,X_{-1} )$ is a positive operator, then $ (A_{-1}+B)_{\vert_{X}} $ generates a positive $ C_0 $-semigroup on $ X $ and
	\begin{displaymath}
		s((A_{-1}+B)_{\vert_{X}})=w_0((A_{-1}+B)_{\vert_{X}})
	\end{displaymath}
\end{theorem}
\begin{proof} 
	To prove our claims we will use \cite[Theorem 2.5]{Arendt}. In fact, let $ M\geq 1 $ and $ w\geq s(A)$ such that $ \Vert T(t)\Vert \leq Me^{wt} $ for all $ t\geq 0 $. Then, for $ \mu >w $, the resolvent of $A$ can be rewrite as 
	\begin{displaymath}
		R(\mu,A)x=\sum_{k=0}^{\infty} e^{-\mu (k\t+\t)}T(k\t) \int_0^{\t} T_{-1}(\t-s) v_{\mu}(s) ds, 
	\end{displaymath} 
	where $ v_\mu=e^{\mu\cdot}x $ for all $ x\in X $. Using the admissibility of $ B $ for $ A $, we then obtain 
	\begin{align*}
		\Vert R(\mu,A)Bx\Vert & \leq\sum_{k=0}^{\infty} e^{-\mu (k\t+\t)}\Vert T(k\t) \Vert \left\Vert\int_0^{\t} T_{-1}(\t-s)B v_\mu(s) ds\right\Vert \\
		&\leq \left( \kappa(\t) \t^{\frac{1}{p}}+\kappa(\t) \t^{\frac{1}{p}}M 
		\frac{e^{(w-\mu)\t}}{1-e^{(w-\mu) \t}}\right)\Vert x\Vert,
	\end{align*} 
	for all $ x\in X $. On the other hand, according to \eqref{S2.R4}, there exists $ \t_0 > 0 $ such that $ \kappa(\t_0) \t^{\frac{1}{p}}_0 $ is sufficiently small.
	Thus, for $ \mu $ sufficiently large,
	\begin{displaymath}
		r(R(\mu,A)B) \leq \Vert R(\mu,A)B\Vert \leq \t_0^{\frac{1}{p}}\kappa(\t_0)<1.
	\end{displaymath}
	Hence, by virtue of \cite[Theorem 3.2]{BJVW} the operator $(A_{-1}+B)_{\vert_{X}}$ is resolvent positive. Furthermore, for $\mu>0$ large enough,
	\begin{align}\label{Resolvent}
		R(\mu,(A_{-1}+B)_{\vert_{X}})=R(\mu,A)+\sum_{n=1}^{\infty} (R(\mu,A_{-1})B)^{n} \ge R(\mu,A).
	\end{align}
	Thus, in view of \eqref{inverse-second} there exist $\mu>s((A_{-1}+B)_{\vert_{X}})$ (sufficiently large) and a constant $c'>0$ such that the operator $R(\mu,(A_{-1}+B)_{\vert_{X}})$ satisfies the inverse estimate \eqref{inverse-second}. Our claims now follow from \cite[Theorem 2.5]{Arendt}. This ends the proof.
\end{proof}

\begin{remark}
	We can rephrase this result by saying that: if $B\in \mathcal{L}_+(X,X_{-1})$ is a positive perturbation of a densely defined resolvent positive operator satisfying the inverse estimate \eqref{S2.5} on a right half-plane, then $B$ is a Desch-Schappacher perturbation for $A$. In addition, according to \eqref{Resolvent}, there exist $\mu_1>s((A_{-1}+B)_{\vert_{X}})$ and $c'>0$ such that 
	\begin{align*}
		\Vert R(\mu_1,(A_{-1}+B)_{\vert_{X}})f \Vert \geq \Vert R(\mu_0, A)f\Vert \ge c' \Vert f \Vert, \qquad \forall f\in X_+.
	\end{align*}
\end{remark}

\section{Admissible positive observation operators}\label{Sec3}
Let $ X,Y $ be Banach lattices. On $ X $ we consider the Cauchy problem
\begin{align*}
	(\mathsf{LOS})  \begin{cases}
		\dot{z}(t) =A z(t),& t> 0,\; z(0)=x,\\
		y(t) = Cz(t),& t> 0.
	\end{cases}
\end{align*}
Here, $ A $ generates a positive C$_0$-semigroup $ \mathbb{T} $ on $ X $ and, $ C $ is positive and maps $ D(A) $ boundedly into $ Y $. We underline that $C$ is not necessarily closed or closeable. Notice that for initial conditions $x\in D(A)$, we have $z(t)=T(t)x$ and then the observation function $y(t)=CT(t)x$ is well-defined for any $t\ge 0$, since the domain $D(A)$ is stable by the semigroup $\mathbb{T}$. However, it is not clear how to define $y(t;x)$ for arbitrary $x\in X$. We say that the system $(\mathsf{LOS}) $ is $L^p$-well-posed (with $p\ge 1$) if the output function $t\mapsto y(t)$ can be extended to a function (denoted by the same symbol) $y\in L^p_{loc}(\R_+,Y)$ such that
\begin{align*}
	\|y(\cdot;x)\|_{L^p([0,\al],Y)}\le \ga \|x\| \qquad (x\in X),
\end{align*}
for any $\al>0$ and some constant $\ga:=\ga(\al)>0$.

We select the following definition.
\begin{definition}\label{S2.D2}
	An operator $ C\in\mathcal{L}(D(A),Y)$ is called $ L^p $-admissible observation operator for $ A $  (with $p\in [1,+\infty)$) if for some (hence for all) $\alpha >0,$ there exists a constant $\gamma:=\gamma(\alpha)>0$ such that
	\begin{align}\label{S2.7}
		\int_{0}^{\alpha} \Vert CT(t)x\Vert^{p}dt\leq \gamma^{p}\Vert x\Vert^{p},\qquad \forall\, x\in D(A). 
	\end{align}
	If $\lim_{\alpha\mapsto 0 }\gamma(\alpha)=0 $, then $C$ is called a zero-class $ L^p $-admissible observation operator.
\end{definition}

In the following lemma we show that the admissibility of a positive observation operator is fully determined on the positive part of $ D(A) $.
\begin{lemma}\label{S2.L6}
	Let $ X , Y $ be Banach lattices, $ \mathbb{T} $ be a positive $ C_0 $-semigroup on $ X $ and $ C\in\mathcal{L}(D(A),Y ) $ be a positive operator. If for some $ \alpha > 0 $ the estimate \eqref{S2.7} holds for any $ 0\leq x\in D(A)$, then $ C $ is an $ L^p $-admissible positive observation operator for $A$.
\end{lemma}

\begin{proof}
	Assume that \eqref{S2.7} holds for any $ 0\leq x\in D(A)$. Let $ x\in D(A) $, then there is $ x_{+},x_{-}\in X_+ $ such that $ x=x_{+}-x_{-} $. Let 
	$
	x_{n,_{\tfrac{+}{}}}:=n\int_{0}^{\frac{1}{n}}T(t)x_{\tfrac{+}{}}dt\in D(A). 
	$
	It is clear that  $ x_{n,_{\tfrac{+}{}}}\to x_{\tfrac{+}{}} $ in $ X $. Further, it follows from
	\begin{displaymath}
		A(x_{n,+}-x_{n,-})=n\int_{0}^{\frac{1}{n}}T(t)A xdt,
	\end{displaymath}
	that $ x_{n,+}-x_{n,-}\to x_{+}-x_{-}  $ in the graph norm of $ A $. Then, we get
	\begin{align*}
		\int_{0}^{\alpha} \Vert CT(t)( x_{n,+}- x_{n,-})\Vert^{p}dt
		\leq \gamma^{p}(\Vert  x_{n,+}\Vert^{p}+\Vert  x_{n,-}\Vert^{p}),
	\end{align*}
	since $0\le x_{n,-},x_{n,+}\in D(A)$. On the other hand, for $ \mu > s(A) $, we have
	\begin{align*}
		&	\Vert CT(t)x-CT(t)(x_{n,+}-x_{n,-})\Vert^{p}\\ 
		&	\leq M^{p}\Vert CR(\mu ,A)\Vert^{p}e^{p w t}\Vert (\mu I_X-A)\big[x-(x_{n,+}-x_{n,-})\big]\Vert^{p}.
	\end{align*}
	It follows that
	\begin{align*}
		\int_{0}^{\alpha}\Vert CT(t)x\Vert^{p} dt
		&\leq \int_{0}^{\alpha}\Vert CT(t)x-CT(t)(x_{n,+}-x_{n,-}) \Vert^{p} dt \\&+ \int_{0}^{\alpha} \Vert CT(t)(x_{n,+}-x_{n,-})\Vert^{p} dt\\
		&\leq  K^{p}(\alpha)\Vert (\mu I_X-A)[(x_{n,+}-x_{n,-})-x]\Vert^{p} \\&+ \gamma^{p}(\alpha)(\Vert  x_{n,+}\Vert^{p}+\Vert  x_{n,-}\Vert^{p}).
	\end{align*}
	Passing to the limit as $ n\mapsto +\infty$, we obtain
	\begin{displaymath}
		\int_{0}^{\alpha} \Vert CT(t)x\Vert^{p} dt \leq \gamma^{p}\Vert x\Vert^{p},
	\end{displaymath}
	for any $ x\in D(A) $ and a constant $ \gamma:= \max\{K(\alpha),\gamma(\alpha)\} >0$.
\end{proof}

\begin{remark}\label{S3.R2}
	To $ C\in\mathcal{L}(D(A),Y)$ we associate an operator $ C_{\Lambda}$, called the \emph{Yosida extension} of $ C $ with respect to $ A $, whose domain, denoted by $ D(C_{\Lambda}) $, consists of all $ x\in X $ for which $\lim_{\mu \mapsto \infty}C\mu  R(\mu ,A)x$ exists ($\mu\in \R$). If $ C $ is positive and $ A $ is resolvent positive, then by the closedness of the positive cone $Y_+$ we have $  C_{\Lambda} $ is also positive.  Then, according to \cite{WO}, for an $ L^p $-admissible positive observation operator $C $, we have $ 0 \leq T(t)x\in D(C_{\Lambda}) $ and
	\begin{align}\label{S2.C}
		(\Psi x)(t)= C_{\Lambda}T(t)x \geq 0,
	\end{align}
	for all $ x\in X_+ $ and for a.e. $ t\geq 0 $. Furthermore, if $ \Psi_{\tau}x=\Psi x $ on $ [0,\tau] $ for any $x \in X$ and $ \tau\geq 0$, then $((\Psi_{\tau})_{\t \geq 0},\mathbb{T})$ define a positive linear observation system on $ X , Y $, see \cite{WO} for more details.
	
	In view of \cite[Theorem 1.10]{CHARALAMBOS}, it is to be noted that an abstract positive linear observation system can be seen as the extension of a family of operators  $ \Psi_\t: X_+\to L^{p}_+(\R_+;Y) $ satisfying
	\begin{align}\label{S3.comp1}
		\begin{array}{lll}
			\Psi_\t (x+y) &=\Psi_\t x
			+\Psi_\t y \\
			\Psi_{\t+t}x &= \Psi_{t}x \qquad\qquad \qquad\quad\,
			{\rm on }\;\; [0,t] \\
			\Psi_{\t+t}x &= [\Psi_{\t}T(t)x](\cdot-t) \qquad {\rm on }\; \; [t,\t+t],
		\end{array}
	\end{align}
	for $\t,t\geq 0$ and $x,y\in X_+$. Note also that, according to \cite[Theorem 3.3]{WO}, for a positive observation system $(\Psi,\mathbb{T})$ there exists a unique positive $L^p$-admissible observation operator $C$. In fact, the positivity of $C$ follows from the fact that
	\begin{displaymath}
		Cx=\lim_{\t\mapsto 0}\frac{1}{\t}\int_0^\t (\Psi x)(t)dt,
	\end{displaymath}
	for all $0\leq x\in D(A)$ (see \cite[Formula 3.4]{WO}) and the closedness of $Y_+$.
\end{remark}

In view of the above lemma, the following result can be obtained by a slight modification of the proof of Lemma 2.1 in Voigt \cite{Voigt}.
\begin{lemma}\label{S2.P3}
	Let $ X , Y $ be Banach lattices such that $ Y $ is a real AL-space. Assume that $ \mathbb{T} $ is a positive $ C_0 $-semigroup on $X$ and $ C\in\mathcal{L}(D(A),Y ) $ is positive. Then, $ C $ is an $L^1$-admissible positive observation operator for $ A $.
\end{lemma}
\begin{remark}
	Note that, according to \cite[Theorem 3.3]{WO}, if $(\Psi,\mathbb{T})$ is a positive observation system, then there exists a unique positive $L^p$-admissible observation operator $C$. In fact, the positivity of $C$ follows from the fact that
	\begin{displaymath}
		Cx=\lim_{\t\mapsto 0}\frac{1}{\t}\int_0^\t (\Psi x)(t)dt,
	\end{displaymath}
	for all $0\leq x\in D(A)$ (see \cite[Formula 3.4]{WO}) and the closedness of $Y_+$.
\end{remark}
Now, we provide a sufficient condition for admissibility of positive observation operators for positive semigroups.
\begin{proposition}\label{S2.P4} 
	Let $ X,Y $ be Banach lattices and let $ A$ generate a positive C$_0$-semigroup on $X$. Furthermore, we assume that $ A^* $ is densely defined and there exist $\mu_0>s(A^*) $ and $c>0$ such that
	\begin{align*}
		\Vert R(\mu_0,A^*)x^* \Vert\geq c \Vert x^* \Vert,\qquad (x^*\in X'_+),
	\end{align*}
	If $C\in\mathcal{L}(D(A),Y) $ is a positive operator, then $ C $ is an $ L^p $-admissible positive observation operator for $ A $.
\end{proposition}
\begin{proof}
	First, we note that $X'$ the topological dual of $X$ is also a Banach lattice and $R(\mu,A)^*=R(\mu,A^*)$. Then, according to \cite[Theorem 2.5]{Arendt}, the operator $ A^* $ generates a positive $ C_0 $-semigroup on $ X' $. Now, let $ \alpha>0 $ and $ 0\leq x\in D(A)$, then using Theorem \ref{S2.T1}  we obtain
	\begin{align*}
		\int_0^\alpha \Vert CT(t)x\Vert^{p} dt&=\sup_{u\ge 0,\Vert u\Vert\leq 1}
		\int_{0}^{\alpha} \langle CT(t)x,u(t)\rangle dt\\
		&=\sup_{u\ge 0,\Vert u\Vert\leq 1}
		\langle x,\int_{0}^{\alpha} T^*(t)C^*u(t) dt\rangle\\
		&\le \left\Vert C^*\right\Vert_{\mathfrak{B}_q} \Vert x\Vert\\
		&\le \tfrac{1}{c}\Vert  R(\mu_0,A_{-1}^*)C^*\Vert M \alpha^{\frac{1}{p}}e^{w\alpha} \Vert x\Vert,
	\end{align*}
	where $\Vert u\Vert:=\Vert u\Vert_{L^p([0,\alpha];U)}$. Therefore, according to Lemma \ref{S2.L6}, $ C $ is an L$^p$-admissible positive observation operator with $\gamma(\alpha)=\tfrac{1}{c}\Vert  R(\mu_0,A_{-1}^*)C^*\Vert M \alpha^{\frac{1}{p}}e^{w\alpha}$. 
\end{proof}
\begin{remark}\label{S3.R3}
	In the presented form, the above proposition asserts that if the adjoint operator $ A^* $ is densely defined and that satisfies the inverse estimate \eqref{S2.5}, then $ C \in\mathcal{L}_+(D(A),Y ) $ is a zero-class $ L^p $-admissible positive observation operator for $ A $ for any $p\in [1,+\infty)$ as $\lim_{\t\mapsto 0 } \gamma(\alpha)=0$.
\end{remark}

\section{Positivity of input-output system}\label{Sec4}
Let $ X,U,Y $ be Banach lattices. Let us consider the input-output linear system 
\begin{align*}
	(\Sigma) \begin{cases}
		\dot{z}(t) = A_{-1} z(t)+Bu(t),& t> 0,\; z(0)=x,\\
		y(t) = \Gamma z(t)+Du(t),& t>  0.
	\end{cases}
\end{align*}
Here $ A $ generates a strongly continuous positive semigroup $ \mathbb{T} $ on $ X $, $ B\in\mathcal{L}(U,X_{-1})$ is a (positive) control operator, and $ D\in\mathcal{L}(U,Y)$ is a (positive) operator, and $ \Gamma\in\mathcal{L}(Z,Y)$ is a (positive) observation operator with
\begin{align*}
	Z=X_1+R(\mu,A)BU, \qquad \mu>s(A),
\end{align*}
which is an order Banach space with the norm 
\begin{align*}
	\Vert z\Vert_{Z}^2=\inf\big\{\Vert x\Vert^2+\Vert u\Vert^2: x\in X_1,\,u\in U,\;z=x+R(\mu,A) Bu \big\}.
\end{align*}
Note that $Z$ is independent of the choice of $\mu$ since different choices of $\mu$ lead to equivalent norms on $Z$ due the resolvent equation. Note also that $X_1\subset Z\subset X$ with continuous embedding. 

We recall the following definition (see, e.g., \cite{Salam,Staf,WR}).
\begin{definition}\label{Definition-Lp-well-posed}
	Let $ X,U,Y $ be Banach spaces. We say that $(\Sigma) $ is an $L^p$-well-posed system (with $p\ge 1$) if for every $ \tau > 0 $ there exists a constant $ c_\tau>0 $ (independent of $u$ and the initial state $x$) such that the following inequality holds for all solutions of $(\Sigma) $:
	\begin{align}\label{Sigma}
		\Vert z(\tau)\Vert^{p}_{X}+\Vert y(\cdot)\Vert^{p}_{ L^{p}([0,\tau];U)}\leq c_\tau\left( \Vert x\Vert^{p}_{X}+\Vert u(\cdot)\Vert^{p}_{ L^{p}([0,\tau];U)}\right).
	\end{align}
\end{definition}
According to Section \ref{Sec2}, the fact that $B$ is an $L^p$-admissible positive control operator for $A$ implies that the state of $(\Sigma)$ is a continuous $ X $-valued function of $ \t $ and satisfies
\begin{align*}
	0\le z(\t)=T(\t)x+\Phi_{\t}u,\qquad \forall\,(x,u) \in X_+\times L^{p}_{+}(\R_+;U) .
\end{align*}
Then, the observation function $ y(\cdot) $ is not well defined as $\Gamma$ is unbounded and $\Phi_\t u\in X$ for all $\tau>0$. To avoid the aforementioned problem, let us consider the following space 
\begin{align*}
	W^{1,p}_{0}([0,\t];U):=\left\{u\in W^{1,p}([0,\t];U):u(0)=0\right\},
\end{align*}
for a fixed $\tau>0$. Assuming $0\in \rho(A)$ (without loss of generality) and using an integration by parts argument, one can prove that $\Phi_\t u\in Z $ for any $u\in W^{1,p}_{0}([0,\tau],U)$. Therefore, for $u$ as before, it makes sense to define the operator
\begin{align*}
	(\mathbb{F}u)(t):=\Gamma\Phi_t u +Du(t), \qquad t\in [0,\t],
\end{align*}
for all $u\in W^{1,p}_{0}([0,\tau],U)$.

We select the following definition.
\begin{definition}\label{S4.D1}
	Let $ X,U,Y $ be Banach spaces. Let $A$ generate a $ C_0 $-semigroup $ \mathbb{T} $ on $ X $, $B$ be an $L^p$-admissible control operator for $ A $, and $C:=\Gamma_{\vert D(A)}$ be an $L^p$-admissible observation operator for $ A $. We say that $(A,B,C) $ is an $ L^p $-well-posed triplet on $(X,U,Y)$ (with $p\ge 1$) if for every $ \tau >0 $ there exists $ \beta:=\beta(\tau) >0$ such that
	\begin{align}\label{Condition1}
		\Vert \mathbb{F} u\Vert_{L^{p}([0,\tau];Y)}\leq \beta \Vert u\Vert_{L^{p}([0,\tau];U)},
	\end{align}
	for all $u\in W^{1,p}_{0}([0,\t];U)$. Here, $\mathbb{F}$ is called the input-output map of $(A,B,C)$. 
\end{definition}
In view of the above definition we obtain the following result.
\begin{lemma}\label{ABC-sigma}
	Let $ X,U,Y $ be Banach spaces and let $(A,B,C) $ be an $ L^p $-well-posed triplet on $(X,U,Y)$ (with $p\in [1,+\infty)$). Then, the input-output system $(\Sigma)$ is $L^p$-well-posed.
\end{lemma}
\begin{proof}
	According to Definition \ref{Definition-Lp-well-posed}, we have to show that the state $z(\cdot)$ and the output $y(\cdot)$ of $(\Sigma)$ satisfy the estimate \eqref{Sigma}. In fact, let $ M\geq 1 $ and $ w\geq w_0(A)$ such that $ \Vert T(t)\Vert \leq Me^{wt} $ for all $ t\geq 0 $. It follows from the $L^p$-admissibility of $B$ that $(\Sigma)$ has a unique mild solution $z(\t)=T(\t)x+\Phi_{\t}u$ for all $\tau\ge 0$. Then, for every $\tau\ge 0$,
	\begin{align*}
		\|z(\t)\|^p_X\le M_1\left(\|x\|_X^p+\|u\|_{L^p([0,\t];\partial X)}^p\right), 
	\end{align*}  
	for all $(x,u)\in X\times L^p([0,\tau];U)$, where $M_1:=M_1(\tau)>0$ dependent of $\kappa,M,w$. On the other hand, the output function
	\begin{align*}
		y(t)=\Gamma z(t)+D u(t)&=CT(t)x+\Gamma \Phi_{\t}u+Du(t)\\&:=(\Psi x)(t)+(\mathbb{F}u)(t)
	\end{align*}
	is well-defined for all $t\ge 0$ and $(x,u)\in D(A)\times W^{1,p}_0([0,\t],\partial X)$, where $C:=\Gamma_{\vert D(A)}$. The fact that $C$ is an $L^p$-admissible observation operator for $A$ implies that $\Psi$ is extendable to $\Psi\in \mathcal{L}(X,L^{p}([0,\t];\partial X))$ for all $\t\ge 0$ such that $(\Psi x)(t)=C_{\Lambda}T(t)x$ for all $x\in X$ and a.e. $t\in [0,\t]$. The estimation \eqref{Condition1} yields that the input-output map $ \mathbb{F}$ is uniquely extendable to a bounded linear operator (still denoted by $ \mathbb{F}$) from $L^{p}([0,\t];\partial X)$ to $L^{p}([0,\t];\partial X)$ for each $\t\ge 0$. Thus, for every $\tau\ge 0$, 
	$$\Vert y(\cdot)\Vert^{p}_{ L^{p}([0,\t];\partial X)}\leq M_2\left( \Vert x\Vert^{p}_{X}+\Vert u(\cdot)\Vert^{p}_{ L^{p}([0,\t];\partial X)}\right),$$ 
	for all $(x,u)\in X\times L^p([0,\tau];U)$, where the constant $M_2:=M_2(\tau)>0$ dependent of $\gamma,\beta,M,w$.
	This ends the proof.
\end{proof}

\begin{remark}\label{S4.R1}
	It follows from Lemma \ref{ABC-sigma} that the output function $ y $ of $(\Sigma)$ is extended to $ y\in L_{loc}^{p}(\R_+;Y) $ and satisfies
	\begin{align}\label{S2.9}
		y=\Psi x+\mathbb{F} u,\qquad\qquad \forall \,(x,u)\in X\times L_{loc}^{p}(\R_+;U),
	\end{align}
	Note that $\mathbb{F}$ commutes with the right shift operator on $X$. Then, for every $\t\geq 0$, the restrictions
	$
	\mathbb{F}_{\tau}=\mathbb{F}_{\vert_{[0,\tau]}},
	$
	are well defined. These operators are called the input-output maps of the triplet $(A,B,C)$. Further, there exists $ \alpha> w_0(A) $ and a unique bounded analytic function $ \mathbf{H}:\mathbb{C}_{\alpha}\to \calL(U,Y) $ such that 
	\begin{align}\label{S2.10}
		\hat{y}(\mu)&=CR(\mu,A)x+\mathbf{H}(\mu) \hat{u}(\mu)
	\end{align}
	for any $(x,u)\in X\times L^{p}_{loc}(\R_+;U)$ and $\Re e\,\mu> \alpha $. Here $ \mathbf{H} $ denotes the transfer function of $(A,B,C) $, see e.g., \cite[Chap. 4]{Staf} for more details. If, in addition,
	\begin{align}\label{S4.d}
		[D u=]\quad	\lim_{\mathbb{R}\ni \mu\mapsto +\infty }\mathbf{H}(\mu)u,
	\end{align}
	exists for all $u\in U$ in the strong (resp., weak) topology of $Y$, then $(A,B,C)$ is called an $ L^p $-well-posed strongly regular (resp., weakly regular) triplet on $(X,U,Y)$. If $(A,B,C)$ is an $ L^p $-well-posed strongly regular triplet, then the state $z$ and the output $y$ of $(\Sigma)$ satisfy
	\begin{subequations}
		\begin{align}
			z(t)&=T(t)x+\Phi_{t}u\in D(C_{\Lambda}) \label{eq:aa}\\
			y(t)&=C_{\Lambda}z(t)+Du(t) \label{eq:bb},
		\end{align}
	\end{subequations}
	for all $ (x,u)\in X\times L_{loc}^{p}(\R_+;U) $ and a.e. $ t \ge 0 $, cf. \cite[Theorem 5.6.5]{Staf}. 
\end{remark}

We can then select the following two definitions, which are the counterparts of \cite[Definition 1]{FR} and \cite[Definition 2]{FR} in the infinite-dimensional setting.
\begin{definition}\label{S4.D2}
	Let $(A,B,C) $ be an $ L^p $-well-posed triplet, on $X,U,Y$ (Banach lattices), in the sense of Definition \ref{S4.D1}. We say that the system $(\Sigma) $ is positive (or internally positive), if for every positive initial state and for every positive input function the state and the output function of $(\Sigma) $ remain positive. 
\end{definition}

\begin{definition}\label{S4.D3}
	Let $(A,B,C) $ be an $ L^p $-well-posed triplet, on $X,U,Y$ (Banach lattices), in the sense of Definition \ref{S4.D1}. We say that the system $(\Sigma) $ is externally positive, if the output function corresponding to the zero initial state is positive for every positive input function. 
\end{definition}
Note that the positivity of the input and the output functions in the above two definitions is understood in the sense of $L^p$-functions. 

In the following result we derive a simple characterization of externally positive linear system. 
\begin{proposition}\label{S4.P1}
	Assume that $(A,B,C) $ is an $ L^p $-well-posed triplet on  $ X,U,Y $ (Banach lattices) with $p\in [1,+\infty)$. Then the following assertions are equivalent:
	\begin{itemize}
		\item[$(\emph{i})$] The system $(\Sigma) $ is externally positive.
		\item[$(\emph{ii})$] The input-output map $ \mathbb{F}$ is positive.
		\item[$(\emph{iii})$] The input-output maps $ \mathbb{F}_t$ are positive for every $t\geq 0$.
	\end{itemize}
\end{proposition}
\begin{proof}
	The equivalence of (\emph{i}) and (\emph{ii}) follows directly from \eqref{S2.9} (with $x=0$). On the other hand, the fact that the input-output maps $ \mathbb{F}_t$ are defined by truncating the input-output map $ \mathbb{F}$ to $[0,t]$ yields the equivalence of (\emph{ii}) and (\emph{iii}).
\end{proof}

On the other hand, the internal positivity of $(\Sigma) $ is characterized as follows. 
\begin{proposition}\label{S4.P2}
	Let $(A,B,C) $ be an $ L^p $-well-posed triplet on $ X,U,Y $ (Banach lattices) with $p\in [1,+\infty)$. Then the system $(\Sigma) $ is internally positive if and only if  $A$ generates a positive $C_0$-semigroup $\mathbb{T}$ on $X$, $C$ is positive, $B$ is positive, and $\mathbb{F}$ is positive. If, in addition, $(A,B,C) $ is an $ L^p $-well-posed strongly regular triplet, then $(\Sigma) $ is internally positive if and only if $D,C,B,\mathbb{T}$ are positive.
\end{proposition}	
\begin{proof}
	Let $(A,B,C) $ be an $ L^p $-well-posed triplet. Assume that $(\Sigma) $ is internally positive. Then, taking $u$ identical zero in \eqref{eq:aa} yields that $z(t)=T(t)x\in X_+$ for any $x\in X_+$ and $t\geq 0$, i.e., $\mathbb{T}$ is positive. Letting $x=0$ in \eqref{eq:aa} for arbitrary vector $v\in U_+$ and $t\geq 0$, we obtain
	$
	\tfrac{1}{t}\Phi_tv= \tfrac{1}{t}z(t) \in X_+,
	$
	which tends to $Bv$ as $t\mapsto 0$. The closedness of $X_{-1,+}$ implies $Bu\in X_{-1,+}$ (recall that $ X_{+}=X\cap  X_{-1,+}$). Moreover, we have $Cx=Cz(0)=y(0)\geq 0 $ for any $x\in X_+$ and $u\equiv 0$ and hence $C$ is positive. Finally, the fact that $0\leq y=\mathbb{F}u$ implies that $\mathbb{F}$ is positive.
	
	Conversely, let $\mathbb{F},C,B$ be positive operators and $\mathbb{T}$ be a positive $C_0$-semigroup on $X$. By virtue of Lemma \ref{S2.LB} the input maps $\Phi_t$ are positive for every $t\geq 0$. This yields that $z(t)$ is positive for every $t\geq 0$. Moreover, it follows from \eqref{S2.9} that $y(t)\geq 0$ for almost every $t\geq 0$. Thus, in view of Definition \ref{S4.D2}, the system $(\Sigma) $ is internally positive. The last statement now follows from the fact that
	\begin{align}\label{Fa}
		(\mathbb{F} u)(t)=C_{\Lambda}\Phi_t u+ Du(t),
	\end{align}
	for any $u\in L_{loc}^{p}(\R_+;U) $ and a.e. $t\geq 0$, since $(A,B,C) $ is a regular triplet. 
\end{proof}

The following theorem shows that an $ L^p $-well-posed positive triplet $(A,B,C)$ is fully determined by its behavior on the positive cone.
\begin{proposition}\label{S4.T1}
	Let $ X,U,Y $ be Banach Lattices, let $\Gamma,D$ be positive and $p\in [1,+\infty)$. Assume that $A$ generates a positive $ C_0 $-semigroup $ \mathbb{T} $ on $ X $, $B$ is an $L^p$-admissible positive control operator for $ A $ and $C$ is an $L^p$-admissible positive observation operator for $ A $. If, in addition, for every $ \tau >0 $ there exists $ \beta:=\beta(\tau)>0 $ such that
	\begin{align}\label{S4.Fu}
		\Vert \mathbb{F} u\Vert_{L^{p}([0,\tau];Y)}\leq \beta \Vert u\Vert_{L^{p}([0,\tau];U)},
	\end{align}
	for all $0\leq u\in W^{1,p}_{0}([0,\t];U)$.	Then $(A,B,C)$ is a positive $L^p$-well-posed triplet on $X,U,Y$. In particular, the system $(\Sigma) $ is internally positive.
\end{proposition}
\begin{proof}
	Let $\t>0$. Since the positive cone of $W^{1,p}_{0}([0,\t];U)$ is dense in the positive cone of $L^{p}([0,\tau];U)$, then the estimate \eqref{S4.Fu} is valid for all $0\leq u \in L^{p}([0,\tau];U)$. Moreover, the fact that 
	\begin{displaymath} 
		(\mathbb{F} u)(t)=\Gamma\Phi_t u+Du(t),\qquad t\in[0,\t]
	\end{displaymath}
	is positive implies that $ \mathbb{F}: L^{p}_+([0,\tau];U)\to L^{p}_+([0,\tau];Y)$. Thus, according to \cite[Theorem 1.10]{CHARALAMBOS}, the operator $ \mathbb{F}$ is uniquely extendable to a positive bounded linear operator $ \mathbb{F}\in \calL (L^{p}([0,\t];U),L^{p}([0,\t];Y)) $ for each $\tau\ge 0$. Therefore, in view of Definition \ref{S4.D1} we get that $(A,B,C)$ is a positive $ L^p $-well-posed triplet. The last claim follows from Proposition \ref{S4.P2}. This completes the proof.
\end{proof}

The next result present an interesting properties of $L^p$-well-posed positive systems.
\begin{theorem}\label{S4.P4}
	Let $ X , U,Y$ be Banach lattices and let $(A,B,C)$ be a positive $ L^p $-well-posed triplet on $ X,U,Y $. Then the following assertions are equivalent:
	\begin{itemize}
		\item[$(\emph{i})$]  $(A,B,C) $ is a positive $ L^p $-well-posed strongly regular triplet.
		\item[$(\emph{ii})$] $(A,B,C) $ is a positive $ L^p $-well-posed weakly regular triplet.
	\end{itemize}
	In this case, the feedthrough operator $D$ is positive.
\end{theorem}
\begin{proof}
	The implication (\emph{i})$\implies$(\emph{ii}) is straightforward. If (\emph{ii}) holds then the transfer function $\mathbf{H}$ of $(A,B,C)$ converge weakly to $D$ and satisfies
	\begin{align*}
		\frac{\mathbf{H}(\mu)-\mathbf{H}(\lambda)}{\mu -\lambda}= -CR(\mu,A)R(\lambda,A)B,
	\end{align*}
	for any $\mu,\lambda>s(A)$ with $\mu\neq\lambda$, cf. \cite[Equation (4.13)]{WR}. Thus, 
	\begin{align}\label{transfequation}
		\mathbf{H}(\mu)v\le \mathbf{H}(\lambda)v, \qquad \forall\, \la \le \mu,\, v\in U_+.
	\end{align}
	Let $v\in U_+$ and let $(\mu_n)_{n\in \mathbb{N}}$ be a increasing sequence in $(s(A),+\infty)$ tending to $+\infty$. It follows from \eqref{transfequation} that $(\mathbf{H}(\mu_n)v)_{n\in \mathbb{N}}$ is monotone decreasing and so, by \cite[Proposition 10.9]{BFR}, we have strong convergence. The last claim follows from the closedness of $Y_+$.
\end{proof}

In particular, we obtain the following result. 
\begin{corollary}\label{S4.C2}
	Every positive $ L^1 $-well-posed triplet on Banach lattices is strongly regular.
\end{corollary}
\begin{proof}
	The proof follows from \cite[Theorem 5.6.6]{Staf} and Theorem \ref{S4.P4}.
\end{proof}

Next, we confine our attention to closed-loop positive systems. Indeed, let us assume that $(A,B,C)$ is a positive $ L^p $-well-posed triplet on $X,U,Y$ (Banach lattices). Consider the output-feedback law 
$
u(t)=Ky(t)+v(t),
$
where $u$ is the input function of $(\Sigma) $, $y$ is the output function of $(\Sigma) $, $v$ is the new input function, and $K\in \calL(Y,U)$ is the feedback operator. Then, according to Remark \ref{S4.R1}, the feedback law $ u=K y+v $ has a sense if only if $ (I_U-K \mathbb{F} )u=K\Psi x+v$ has a unique solution $ u\in L^{p}_{loc}(\R_+;U) $. This is true if $ I_U-K \mathbb{F} $ is invertible. In this case, we say that $K$ is an admissible feedback operator for $ (A,B,C) $. This new system is called the \emph{closed-loop system} corresponding to $(\Sigma) $ and $K$, and it is denoted by  $(\Sigma^K)$, see e.g., \cite[Chap. 7]{Staf} for more details.

Thus we introduce the following definition.
\begin{definition}\label{S4.D4}
	Let $(A,B,C)$ be a positive $ L^p $-well-posed triplet on $X,U,Y$ (Banach lattices) with $p\in [1,+\infty)$. We say that $K\in \calL(Y,U)$ is a positive admissible feedback operator for $ (A,B,C) $, if $ I_U-K\mathbb{F}$ has a positive inverse in $\calL(L^{p}([0,\tau];U)) $ for some $\tau>0$.
\end{definition}
We have the following characterization of positive admissible feedback operators.
\begin{lemma}
	Let $X,U,Y$ be Banach lattices and let $(A,B,C)$ be a positive $ L^p $-well-posed triplet on $X,U,Y$ (with $p\in [1,+\infty)$). Then, $K$ is a positive admissible feedback operator for $ (A,B,C) $ if and only if $r(K \mathbb{F})<1 $.
\end{lemma}
\begin{proof}
	Necessity, let $K\in \calL(Y,U)$ be a positive admissible feedback operator for $ (A,B,C) $. Then, according to Definition \ref{S4.D4}, $ (I_U-K\mathbb{F})^{-1}$ exists in $\calL(L^{p}([0,\t];U)) $ for some $\t>0$ and $ (I_U-K\mathbb{F})^{-1}\geq 0$. It follows that $r(K \mathbb{F}) \leq 1$. On the other hand, the positivity of $K\mathbb{F}$ implies $r(K \mathbb{F})\in \sigma(K\mathbb{F}) $ by Pringsheim's theorem (cf. \cite[Appendix 2.2]{Schaf1}), and therefore $r(K \mathbb{F})<1 $.
	
	Sufficiency follows directly from the fact $(I_U-K\mathbb{F})^{-1}$ can be expanded in Neumann series
	$$
	(I_U-K\mathbb{F})^{-1}=\sum_{n=0}^\infty (K\mathbb{F})^n \geq I_U \geq 0.
	$$
\end{proof}
\begin{remark}
	It is easy to see that the above proposition remain valid if we replace $r(K \mathbb{F})<1 $ by $r(\mathbb{F} K )<1 $. Indeed, we have $(I_Y+\mathbb{F}(I_U-K\mathbb{F})^{-1}K)(I_Y-\mathbb{F}K)=I_Y$ and $ (I_Y-\mathbb{F}K)(I_Y+\mathbb{F}(I_U-K\mathbb{F})^{-1}K)=I_Y $, hence $ (I_Y-\mathbb{F}K)^{-1}= I_Y+\mathbb{F}(I_U-K\mathbb{F})^{-1}K \geq 0$.
\end{remark}

\begin{lemma}\label{S4.L1}
	Let $X,U,Y$ be Banach lattices such that $Y$ has an order continuous norm. Let $(A,B,C)$ be a positive $ L^p $-well-posed strongly regular triplet on $X,U,Y$ and $K\in \calL(Y,U)$ be a positive operator. If $ r(K\mathbf{H}(\mu_0))< 1$ for some $\mu_0> s(A)$, then we obtain   
	\begin{displaymath}
		r(KD)< 1\qquad {\rm and}\qquad r(DK)<1.
	\end{displaymath}
\end{lemma}
\begin{proof}
	According to the proof of Theorem \ref{S4.P4}, we have $(\mathbf{H}(\mu))_{\mu> s(A)}$ is a decreasing subset of $Y$. The fact that $(A,B,C)$ is a positive $ L^p $-well-posed strongly regular triplet further yields that
	\begin{displaymath}
		Du=\inf\{\mathbf{H}(\mu)u:\; \mu>s(A)\}
	\end{displaymath}
	for all $u\in U$, since $Y$ is Banach lattice. Then we have
	$
	KDu \leq K \mathbf{H}(\mu)u,
	$
	for all $u\in U$ and $\mu>s(A)$ (as $K$ is positive), and consequently
	\begin{align}\label{S4.1}
		\sum_{n=0}^N (KD)^n u \leq \sum_{n=0}^N (K \mathbf{H}(\mu_0))^n u \leq (I_U-K \mathbf{H}(\mu_0))^{-1}u, \qquad\forall\, u\in U, 
	\end{align}
	since $ r(K\mathbf{H}(\mu_0))< 1$. The positivity of $KD$ implies that $\sum_{n=0}^N (KD)^n u$ is monotone increasing and so
	\begin{displaymath}
		\sum_{n=0}^\infty (KD)^n u < +\infty, 
	\end{displaymath}
	for all $u\in U$ as we assumed that $Y$ has order continuous norm. Moreover, from the Banach-Steinhaus
	theorem we obtain the boundedness of $\{(KD)^n\}_{n\geq 0}$, so that the series converges also in the uniform operator topology. Therefore, we obtain $ r(KD)< 1$.
\end{proof}

In order to state rigorously the main result of this section, let us consider the following space
\begin{align*}
	\mathscr{W}_{-1}&:=(\la I_X-\tilde{A}_{-1})D(C_\Lambda), \quad \Vert x\Vert_{\mathscr{W}_{-1}}= \Vert R(\la,\tilde{A}_{-1}) x\Vert_{D(C_\Lambda)},
\end{align*}
for $\la\in \rho(A)$, where $\tilde{A}_{-1}$ is the restriction of $A_{-1}$ to $D(C_\Lambda)$ and
\begin{align*}
	\Vert x\Vert_{D(C_\Lambda)}=\Vert x\Vert_X+ \sup_{s>1+\omega_0(A)}\Vert s C R(s,A)x \Vert_Y.
\end{align*}
It follows from \cite[Lemma 5.1.3]{Staf} that $(\mathscr{W}_{-1},\Vert \cdot\Vert_{\mathscr{W}_{-1}})$ is a Banach space independent of the choice of $\la $. In addition, $R(\la,\tilde{A}_{-1})\in \calL(\mathscr{W}_{-1},D(C_\Lambda))$ and $X\subset \mathscr{W}_{-1}\subset X_{-1}$ with continuous embeddings. 

We can now state and prove the main result of this section.
\begin{theorem}\label{S4.T3}
	Let $X,U,Y$ be Banach lattices such that $Y$ has an order continuous norm. Let $(A,B,C)$ be a positive $ L^p $-well-posed strongly regular triplet on $(X,U,Y)$ with $p\in [1,+\infty)$, and let $D$ be the corresponding feedthrough operator. Let $K\in\calL(Y,U)$ be a positive admissible feedback operator for $ (A,B,C) $. Furthermore, let $(\Sigma^K)$ be the closed-loop system of $(\Sigma)$ with semigroup generator $A^K$, control operator $B^K$, observation operator $C^K$ and feedthrough operator $D^K$. Denote by $\mathscr{W}_{-1}^K$ the analogue of $\mathscr{W}_{-1}$ for $\tilde{A}_{-1}^K$ (i.e., $\mathscr{W}_{-1}^K=(\la I_X-\tilde{A}_{-1}^K)D(C_\Lambda)$ for some $\la\in \rho(A^K)$). If there exists $\mu_0> s(A)$ such that $ r(K\mathbf{H}(\mu_0))< 1$, then the closed-loop triplet $(A^K,B^K,C^K)$ is an $ L^p $-well-posed strongly regular on $(X,U,Y)$ with the extended observation operator $C_{\Lambda}^K\in \calL(D(C_{\Lambda}), U)$ and the corresponding feedthrough operator $D^K$ given by
	\begin{align}\label{S4.2}
		C_{\Lambda}^K= (I_Y-DK)^{-1}C_{\Lambda}, & &
		D^K =D(I_U-KD)^{-1},
	\end{align}
	and the operator $A^K$ is given by 
	\begin{align}\label{S4.3}
		\begin{array}{llll}
			A^K&=\tilde{A}_{-1}+BK(I_Y-DK)^{-1}C_{\Lambda}\\
			D(A^K)&=\{x\in D(C_{\Lambda}):\;\;  (\tilde{A}_{-1}+BK(I_Y-DK)^{-1}C_{\Lambda})x \in X     \}.
		\end{array}
	\end{align}
	Moreover, $\mathscr{W}_{-1}^K=\mathscr{W}_{-1}$ and $B^K,B\in  \calL(U,\mathscr{W}_{-1})$ with
	\begin{align}\label{S4.4}
		B^K= B(I_U-KD)^{-1}.
	\end{align}
	If, in addition, there exists  $\mu_1>s(A)$ such that $r(C_{\Lambda}R(\mu_1,\tilde{A}_{-1}) B^KK)<1$, then $(A^K,B^K,C^K)$ is a positive $ L^p $-well-posed strongly regular triplet  .
\end{theorem}

\begin{proof}	
	According to Lemma \ref{S4.L1}, we have $(I_U-KD)^{-1}$ and $(I_Y-DK)^{-1}$ exist and both are positive. Thus, in view of \cite[Theorem 7.5.3]{Staf}-(iii), $(A^K,B^K,C^K)$ is an $ L^p $-well-posed strongly regular triplet on $(X,U,Y)$ with $D(C_{\Lambda}^K)=D(C_{\Lambda})$ and $\mathscr{W}_{-1}^K=\mathscr{W}_{-1}$. Furthermore, the operators $C^K,D^K,A^K$ and $B^K$ are given by \eqref{S4.2},\eqref{S4.3} and \eqref{S4.4}, respectively. Clearly, the positivity of $B^K,C^K,D^K$ follows from their explicit expressions. So, according to Proposition \ref{S4.T1}, it remains to show that $A^K$ generates a positive C$_0$-semigroup. In fact, according to \cite{WF},  
	\begin{align*}
		R(\mu,\tilde{A}^K_{-1})\left(I_{\mathscr{W}_{-1}}-B^K KCR(\mu,A)\right)=R(\mu,\tilde{A}_{-1}),
	\end{align*}
	for any $\mu\in \C$ with $\Re e\mu>\max\{\omega_{0}(A),\omega_{0}(A^K)\}$, where $\tilde{A}^K_{-1}$ is the restriction of $A_{-1}^K$ to $D(C_\Lambda)$. Denote
	\begin{align*}
		Q(\mu)=CR(\mu,A) \in \calL(X,Y), \qquad \forall \mu>s(A).
	\end{align*}
	Notice that $	Q(\mu)x=C_{\Lambda}R(\mu,\tilde{A}_{-1})x $ for $x\in \mathscr{W}_{-1}$. Since the generalized sequence $(Q(\mu))_{\mu>s(A)}$ is positive and monotonically decreasing and the operators $B^K, K$ are positive, then the generalized sequence $(Q(\mu)B^K K)_{\mu>s(A)}$ is also positive and monotonically decreasing. Thus, 
	\begin{align*}
		r(Q(\mu) B^KK)\le r(Q(\mu_1)B^KK)<1, \qquad \forall \mu\ge \mu_1.
	\end{align*}
	Thus, $1\in  \rho((I_Y-Q(\mu) B^K K)^{-1})$ and $0\le (I_Y-Q(\mu) B^K K)^{-1}$. Therefore, 
	\begin{align*}
		\left(I_{\mathscr{W}_{-1}}-B^K KCR(\mu,A)\right)^{-1}=I_{\mathscr{W}_{-1}}+B^K K (I_Y-Q(\mu) B^K K)^{-1} Q(\mu)\ge 0.
	\end{align*}
	Hence, choosing $\mu_1>\max\{\omega_{0}(A),\omega_{0}(A^K)\}$ yields that
	\begin{align*}
		R(\mu,A^K)=R(\mu,A)\left(I_{\mathscr{W}_{-1}}-B^KKCR(\mu,A)\right)^{-1}\ge R(\mu,A)\ge 0, 
	\end{align*}
	for all $\mu\ge \mu_1$. Therefore, $A^K$ has positive resolvent and hence generates a positive C$_0$-semigroup, cf. \cite[Corollary 2.3]{Arendt}. It follows that $(A^K,B^K,C^K)$ is a positive $ L^p $-well-posed strongly regular triplet on $X,U,Y$. This completes the proof.	
\end{proof}

\begin{remark}\label{Remark}
	(a) If $X,U,Y$ are reflexive Banach lattices, then we can drop the assumptions that $Y$ has order continuous norm and $ r(K\mathbf{H}(\mu_0))< 1$ for some $\mu_0> s(A)$. In fact, if $(A,B,C)$ is a positive $ L^p $-well-posed regular triplet on $X,U,Y$, then the dual triplet $(A^*,C^*,B^*)$ is a positive $ L^p $-well-posed weakly regular on $X,Y,U$. Theorem \ref{S4.P4} further yields that $(A^*,C^*,B^*)$ is strongly regular. Now, the claim follows from \cite[Theorem 7.6.1]{Staf}-(iii). 
	
	(b) Note that in the Hilbert lattice setting the assumption $ r(K\mathbf{H}(\mu_0))< 1$ for some $\mu_0> s(A)$ also implies that $K$ is a positive admissible feedback operator for $(A,B,C)$ and hence yields the results of Theorem \ref{S4.T3}.
\end{remark}

We end this section by the following well-posedness result for non-homogeneous boundary value control problems.
\begin{theorem}\label{S4.CB}
	Let $ X,\partial X,U$ be Banach lattices. On $X$ we consider the non-homogeneous boundary value control problem
	\begin{align}\label{non-homogeneous}
		\begin{cases}
			\dot{z}(t)= A_{m}z(t),& t> 0,\\
			z(0)=x\in X_+,\\
			(G-\Gamma) z(t)=Ku(t), & t> 0,
		\end{cases}
	\end{align}
	where $ A_{m}:D(A_{m}) \to X $ is a linear (differential) operator on $ X $, $K\in \calL_+(U,\partial X)$ and $ G,\Gamma: D(A_m)\to \partial X$ are boundary linear positive operators. We assume that:
	\begin{itemize}
		\item[$(\emph{i})$] The restricted operator $ A\subset A_m $ with domain $ D(A)=\ker G $ generates a positive C$_{0} $-semigroup $ \T$ on $ X $.
		\item[$(\emph{ii})$] The boundary operator $ G $ is surjective and $	D_{\mu}:= \left(G_{\vert \ker ( \mu I_X-A_{m} )}\right)^{-1}$ is positive for every $\mu>s(A)$.
		\item[$(\emph{iii})$] The triplet $ (A,B,C) $ is a positive $L^p$-well-posed strongly regular on $ \partial X,X, \partial X $ (with $p\in [1,+\infty)$) with feedthrough zero and $ I_{\partial X} $ as a positive admissible feedback, where
		\begin{displaymath}
			C:=\Gamma_{\vert D(A)}, \qquad\qquad B:=(\mu I_X -A_{-1})D_\mu, \;\;\mu > s(A).
		\end{displaymath}
	\end{itemize}
	Then the operator
	$\mathscr{A}=A_{m}$ with domain $D(\mathscr{A})=\{x\in D(A_m):\: (G-\Gamma) x=0\}$
	generates a positive C$_0$-semigroup on $ X $ given by
	\begin{align}\label{variation}
		\mathscr{T}(t)x=T(t)x+\int_{0}^{t}T_{-1}(t-s)BC_\Lambda\mathscr{T}(s)xds,\qquad t\geq 0,\;\; ({\rm on}\, X).
	\end{align}
	Moreover, the non-homogeneous boundary value control problem \eqref{non-homogeneous} has a unique mild solution satisfying 
	\begin{align}\label{mild-solution}
		0\le z(t)=	\mathscr{T}(t)x+\int_{0}^{t}\mathscr{T}_{-1}(t-s)BK u(s)ds,
	\end{align}
	for all $t\ge 0$ and $(x,u)\in X_+\times L_+^p(\R_+,U)$.
\end{theorem}
\begin{proof}
	First, we rewrite \eqref{non-homogeneous} as the following boundary input-output system
	\begin{align}\label{input-output}
		\begin{cases}
			\dot{z}(t) =A_{m} z(t),& t> 0,\; z(0)=x,\\
			G z(t)-v(t)= Ku(t),& t> 0,\\
			y(t) = \Gamma z(t), & t> 0,
		\end{cases}
	\end{align}
	with the feedback law 
	\begin{align}\label{feedback}
		"\emph{u=y}".
	\end{align}
	According to \cite[Lemmas 1.2 and 1.3]{Gr}, the assumptions (\emph{i})-(\emph{ii}) imply that the domain of $ A_m $ can be decomposed and related to $ A $ as
	$$
	D(A_m)=D(A) \oplus \ker ( \mu I_X- A_{m} ),\qquad \mu >s(A).
	$$ 
	Thus, $\mu D_\mu v=A_m D_\mu v$ for all $v\in \partial X$ and hence 
	\begin{equation}\label{representation}
		A_{m}=A_{-1}+ BG.
	\end{equation}  
	Then, the boundary input-output system \eqref{input-output} can be reformulated as 
	\begin{equation}\label{distributed-parameter-system}
		\left\lbrace
		\begin{array}{lll}
			\dot{z}(t) =A_{-1} z(t)+Bv(t)+BKu(t),& t> 0,\; z(0)=x,\\
			y(t) = \Gamma z(t), & t> 0.
		\end{array}
		\right.
	\end{equation}
	Now, by virtue of Theorem \ref{S4.T3} the assumptions (\emph{iii}) implies that \eqref{distributed-parameter-system} is equivalent to the following open-loop system 
	\begin{equation}\label{closed-loop}
		\left\lbrace
		\begin{array}{lll}
			\dot{z}(t) =(\tilde{A}_{-1}+BC_\Lambda) z(t)+BKu(t),& t> 0,\\
			z(0)=x.&
		\end{array}
		\right.
	\end{equation}
	On the other hand, it follows from the proof of \cite[Theorem 4.1]{HMR} that $\mathscr{A}=(\tilde{A}_{-1}+BC_\Lambda)$. Thus, according to Theorem \ref{S4.T3}, the operator $\mathscr{A}$ generates a C$_0$-semigroup $\mathscr{T}:=(\mathscr{T}(t))_{t\ge 0}$ on $X$ given by \eqref{variation}. 
	
	Next we will claim that $\mathscr{T}$ is positive. In fact, by \eqref{representation}, we have 
	\begin{align*}
		\mathcal{A} x =A_m x =(A_{-1}+BG)x=(A_{-1}+B\Gamma)x,
	\end{align*}
	for all $x\in D(\mathcal{A})$. Thus, for any $\mu >s(A)$, we have 
	\begin{align*}
		\begin{array}{lll}
			(\mu I_X-\mathcal{A} )x &=(\mu I_X- A_{-1}-B\Gamma)x\\
			&=(\mu I_X- A_{-1})(I_X-D_\mu\Gamma)x\\
			&=(\mu I_X- A)(I_X-D_\mu\Gamma)x
		\end{array}       
	\end{align*}
	for all $x\in D(\mathcal{A})$, since $x-D_\mu\Gamma x\in D(A)$. Thus, for any $\mu>s(A)$,
	\begin{align}\label{carac}
		\mu\in \rho(\mathscr{A}) \iff 1\in \rho(D_\mu\Gamma)\iff 1\in \rho(\Gamma D_\mu).
	\end{align}
	On the other hand, it follows from the Hille-Yosida theorem that there exists $w\in \R$ such that $(w,+\infty)\in \rho(\mathscr{A})$, since $\mathscr{A}$ is a generator of a C$_0$-semigroup on $X$. Thus, there exists $\mu_0>0$ (large enough) such that $\mu_0\in \rho(A)\cap \rho(\mathscr{A})$ and hence by \eqref{carac}  $1\in \rho(\Gamma D_{\mu_0})$. For such $\mu_0$, the Pringsheim's theorem (see, e.g., \cite[Appendix 2.2]{Schaf1}) further yields that $r(\Gamma D_{\mu_0})<1$. Therefore,
	\begin{align}\label{resolvent-eqution}
		\begin{array}{lll}
			R(\mu,\mathscr{A})&=(I_X-D_\mu \Gamma)^{-1}R(\mu,A)\\
			&=(I_X+D_\mu(I_{\partial X}-\Gamma D_\mu)^{-1}\Gamma)R(\mu,A),
		\end{array}       
	\end{align}
	for all $\mu\ge \mu_0$,  since that the generalized sequence $(\Gamma D_\mu)_{\mu >s(A)}$ is positive and monotonically decreasing. Thus, $\mathscr{A}$ has positive resolvent and therefore $\mathscr{A}$ is positive, cf. \cite[Corollary 2.3]{Arendt}. Hence, the boundary input-output system \eqref{input-output} with the feedback law \eqref{feedback} is equivalent to the following distributed control system
	\begin{align}\label{distrubeted-control}
		\left\lbrace
		\begin{array}{lll}
			\dot{z}(t) =\tilde{\mathscr{A}}_{-1}z(t)+\mathscr{B}u(t),\qquad t> 0,\\
			z(0) =x, 
		\end{array}\right.
	\end{align}
	where $\tilde{\mathscr{A}}_{-1}$ is the restriction of $\mathscr{A}_{-1}$ to $D(C_\Lambda)$ and $\mathscr{B}:=BK\in \calL_+(U,\mathscr{W}_{-1})$. Therefore, according to Proposition \ref{S2.P1}, the system \eqref{distrubeted-control}, hence \eqref{non-homogeneous}, has a unique positive mild solution given by \eqref{mild-solution}, since by Theorem \ref{S4.T3} $\mathscr{B}$ is an $L^p$-admissible control operator for $\mathscr{A}$. This ends the proof.
\end{proof}
This corollary actually present a generation result for positive perturbations of positive semigroups, namely, unbounded perturbations of the boundary conditions of a generator of positive semigroup.

\section{Application}\label{Sec5}
Let us consider the following system of transport equations on a closed network:
\begin{align*}
	(\Sigma_{\mathsf{TN}})		
	\begin{cases}
		\dfrac{\partial}{\partial t} z_{j}(t,x,v)= v\dfrac{\partial }{\partial x}z_{j}(t,x,v)+q_{j}(x,v)z_{j}(t,x,v),\quad t\geq 0, \;(x,v)\in \Omega_j,\cr  
		z_{j}(0,x,v)= f_{j}(x,v)\ge 0, \qquad\qquad\quad \qquad\qquad\qquad (x,v)\in \Omega_j,\cr
		\imath^{out}_{ij}z_{j}(t,l_j,\cdot)= \mathsf{w}_{ij}\sum_{k\in \calM} \imath^{in}_{ik}\mathbb{J}_k(z_{k})(t,0,\cdot)+\sum_{l\in \calN_c}\mathsf{b}_{il}u_{l}(t,.), \quad t\geq 0,
	\end{cases}
\end{align*}
for $ i\in \{1,\ldots,N\}:=\calN$, $j\in \{1,\ldots,M\}:=\calM$ and $l\in \{1,\ldots,n\}:=\calN_c$ with $\infty>M\ge N\ge n$, where we set $\Omega_j:=[0,l_j]\times [v_{\min},v_{\max}]$ and $M,N,n\in \N$. The corresponding transport equations are defined on the edges of a finite graph $\mathsf{G}$, whose edges are identified with a collection of $M$-intervals $[0,l_j]$ with endpoints "glued" to the graph structure. The connection of such edges being described by the coefficients $\imath^{out}_{ij},\imath^{in}_{ik}\in \{0,1\}$. The flow velocity along the edges is determined by the function $v$, whereas its absorption is determined by the functions $ q_j(\cdot, \cdot) $. The boundary condition determines the propagation of the flow along the various components of the network. The weights $0\le \textsf{w}_{ij}\le 1$  express the proportion of mass being redistributed into the edges and the non-local operators $\mathbb{J}_k$ describes the scattering at the vertices. Moreover, for $i,l\in  \calN\times \calN_c$, the coefficients $\mathsf{b}_{il}$ denotes the entries of the input matrix $ K $ and $ u_l $ denotes the input functions at the vertices. The system $(\Sigma_{\mathsf{TN}})$ is a generalization of the transport network system studied in \cite{EHR1,Rad}.

Next, we are concerned with the existence and uniqueness of a positive mild solution of $(\Sigma_{\mathsf{TN}})$. To this end, we need to fix some notations from graph theory. Here and in the following, we consider a finite connected metric graph $\mathsf{G}=(\mathsf{V},\mathsf{E})$ and a flow on it (the latter is described by $(\Sigma_{\mathsf{TN}})$). The graph $ \mathsf{G} $ is composed by $ N\in \mathbb{N} $ vertices $ \alpha_{1},\,\ldots,\alpha_N $, and by $ M\in \mathbb{N} $ edges $ \mathsf{e}_{1},\,\ldots,\mathsf{e}_M $ which are assumed to be identified with the interval $ [0,l_j] $ with $l_j>0$. Each edge is parameterized contrary to the direction of the flow of material on them, i.e., the material flows from $l_j$ to $0$. The topology of the graph $\mathsf{G}$ is described by the incidence matrix $\mathcal{I}=\mathcal{I}^{out}-\mathcal{I}^{in}$, where $\mathcal{I}^{out}$ and $\mathcal{I}^{in}$ are \emph{the outgoing incidence} and the \emph{incoming incidence matrices} of $\mathsf{G}$ having entries
\begin{align*}
	\imath^{out}_{ij}:=
	\begin{cases}
		1,\quad \text{if  } \begin{tikzpicture}
			\node (P) at (0,0) {$\mathsf{v}_i$};
			\node (T) at (0.7,0.2) {$\mathsf{e}_j$};
			\node (S) at (1.2,0) {};
			\draw[*->,>=latex] (P) to[=1] (S);
		\end{tikzpicture},
		\\
		0,\quad \text{if not},
	\end{cases}\quad
	\imath^{in}_{ij}:=
	\begin{cases}
		1, \quad \text{if  }  \begin{tikzpicture}
			\node (P) at (0,0) {};
			\node (T) at (0.4,0.2) {$\mathsf{e}_j$};
			\node (S) at (1.2,0) {$\mathsf{v}_i$};
			\draw[->*,>=latex] (P) to[=6] (S);
		\end{tikzpicture},
		\\
		0, \quad \text{if not},
	\end{cases}
\end{align*}
respectively. Replacing $1$ by $\textsf{w}_{ij}\geq 0$ in the definition of $\imath^{out}_{ij}$, we obtain the so-called weighted outgoing incidence matrix $\mathcal{I}^{-}_{\mathsf{w}}:=(\mathsf{w}_{ij}\imath^{out}_{ij})$. In this cases, $\mathsf{G}$ is called a weighted graph and its topology is described via the weighted transposed adjacency matrix $\mathbb{A}:=\mathcal{I}^{in}(\mathcal{I}_{w}^{out})^{\top}$ given by, for $i,l\in \calN$,
\begin{align*}
	(\mathbb{A})_{il}:=
	\begin{cases}
		\mathsf{w}_{lj},& \text{if}\; \exists\, \mathsf{e}_j\begin{tikzpicture}
			\node (P) at (2,0) {};
			\node (R) at (2.26,0.2) {$\mathsf{v}_i$};
			\node (T) at (1.3,0.2) {$\mathsf{e}_j$};
			\node (S) at (4,0) {};
			\node (P') at (0.5,0) {};
			\node (P'') at (2.26,0) {};
			\node (T') at (0.5,0.2) {$\mathsf{v}_l$};
			\draw[*->*,>=latex] (P') to[=0] (P'');
		\end{tikzpicture},
		\\
		0, &\text{if not}.
	\end{cases}
\end{align*}

In what follows, we set $$\Omega:= \underset{j\in\calM}{\prod}\Omega_j,\quad l:=(l_j)_{j\in \calM}, \quad l_m=\underset{ j\in \calM}{ \max}\, l_j.$$ Let us consider the Banach lattices $(X,\|\cdot\|_{X})$, $(\partial X,\|\cdot\|_{\partial X})$ and the Banach space $(\mathbb{W},\Vert \cdot \Vert_{\mathbb{W}})$ defined by
\begin{align*}
	X&:=\prod_{j=1}^M L^{1}(\Omega_j),\qquad  \quad\qquad \qquad \qquad\;\Vert f \Vert_{X} :=\sum_{j=1}^{M}\Vert f_j\Vert_{L^{1}(\Omega_j)},\cr
	\partial X &:= L^1([v_{\min},v_{\max}])^N,   \qquad \;\;  \qquad\qquad\Vert g \Vert_{\partial X} :=\sum_{j=1}^{N}\Vert f_j\Vert_{L^{1}([v_{\min},v_{\max}])},\cr
	\mathbb{W}&:=\displaystyle\prod_{j=1}^M W^{1,1}(\Omega_j),\;  \Vert f \Vert_{\mathbb{W}}=\Vert f\Vert_{X}+ \Vert \partial_x f\Vert_{X}.
\end{align*}
Note that $
X_+=\prod_{j=1}^M L^{1}_+(\Omega_j)$. We introduce the following operator
\begin{align}\label{S1.Am}
	A_{m}f&= v\partial_x f+q(\cdot,\cdot)f,\cr
	D(A_{m})&= \big\{ f \in \mathbb{W}: \; \exists\, j\in \calM \;{\rm such \; that }\; \mathsf{i}_{ij}^{out}\neq 0\; {\rm and}\;\\ &f_j(l_j,v)=\mathsf{w}_{ij}\mathsf{i}_{ij}^{out}x_i,\quad {\rm for \; some}\; x\in \C^N \big\},\nonumber
\end{align}
where $q(\cdot,\cdot):={\rm diag}\,(q(\cdot,\cdot))_{j\in \calM}$.  We set $\mathbb{J}={\rm diag}(\mathbb{J}_j)_{j\in \calM}$ with the operators $\mathbb{J}_j$ are given by 
\begin{align}
	\mathbb{J}_j(f_j)(x,v)&=\int_{v_{\min}}^{v_{\max}}\ell_j(x,v,v')f_j(x,v') dv',\qquad (x,v)\in \Omega_j,\; f\in X,
\end{align} 
with $ 0\le \ell_j\in C([0,l_j]; L^{\infty}([v_{\min},v_{\max}]^2))$ for all $j\in \calM$. We also introduce the input space $U:=L^p([v_{\min},v_{\max}])^n$ and the control operator $K $ is given by 
\begin{align*}
	(Ku)_i(t,\cdot) =\sum_{l\in \calN_c}\mathsf{b}_{il}u_l(t,\cdot) ,\qquad \forall\, i\in \calN,\,t\ge 0,
\end{align*}
where $\mathsf{b}_{il}\ge 0$ for all $i,l\in \calN\times \calN_c$.

In order to apply the results of the previous sections, let us introduce the following assumptions:
\begin{itemize}
	\item[{\bf (A1)}] $0<v_{\min}\le v\le v_{\max}$, $ q_j(\cdot, \cdot)\in L^{\infty}(\Omega_j) $ and $\kappa\le q_j$ for all $j\in \calM$ and for some $-\infty<\kappa<0$.
	\item[{\bf (A2)}] Each vertex has at least one outgoing edge.
	\item[{\bf (A3)}] The weights $\textsf{w}_{ij}$ satisfies 	$\sum_{j\in \calM}\textsf{w}_{ij}=1$ for all $ i\in \calN$.
\end{itemize}
In the above, {\bf (A1)} specifies the transport process along each edges, {\bf (A2)} is equivalent to the statement that $\mathcal{I}^{out}$ is surjective, while {\bf (A3)} implies that the boundary condition exhibits standard \emph{Kirchhoff conditions} and the matrix $\mathbb{A}$ is column stochastic. 
\begin{proposition}\label{final}
	Let the assumptions {\bf (A1)}-{\bf (A3)} be satisfied. Then the transport network system $(\Sigma_{\mathsf{TN}})$ has a unique positive mild solution.
\end{proposition}
\begin{proof}
	To prove our claim we will use Theorem \ref{S4.CB}. To this end, we first introduce the operators $ G,\Gamma: \mathbb{W}\to \partial X$ defined by 
	\begin{align}\label{Operators-boundary}
		(G f)_j&:=	\sum_{i\in \calN}\mathsf{i}_{ij}^{out} f_j(l_j,v),\qquad
		\Gamma f :=\mathcal{I}^{in}(\mathbb{J}f)(0,v),
	\end{align}
	for $v\in [v_{\min},v_{\max}]$ and $f\in \mathbb{W}$. If we set  $z(t)=(z_j(t,\cdot,\cdot))_{j\in \calM}$ and $u(t)=(u_l(t,\cdot))_{l\in \calN_c}$, then the non-homogeneous boundary value control problem
	\begin{align*}
		\begin{cases}
			\dot{z}(t)= A_{m}z(t),& t> 0,\\
			z(0)=f\ge 0,\\
			(G-\Gamma) z(t)=Ku(t), & t> 0,
		\end{cases}
	\end{align*}
	with $f=(f_j(\cdot,\cdot))_{j\in \calM}$ is an abstract version of $(\Sigma_{\mathsf{TN}})$. By Theorem \ref{S4.CB} this problem is well-posed if the homogeneous problem is well-posed, namely, if the operator 
	\begin{align}
		\mathscr{A}=v\partial_x f+q(\cdot,\cdot)f, \qquad D(\mathscr{A})=\{x\in D(A_m):\: (G-\Gamma) f=0\},
	\end{align}
	generates a positive C$_0$-semigroup on $X$. Indeed, assumption {\bf (A1)} yields that the restriction $(\mathscr{A}, D(\mathscr{A}))$ to $\ker G$, denoted by $A$, generates a positive C$_0$-semigroup $(T(t))_{t\geqslant 0} $ on $X$ given by
	\begin{align}\label{Semigroup}
		(T(t)f)_j(x,v)=
		\begin{cases}
			e^{\int_0^{t}q_j(x+v\sigma,v)d\sigma}f_j(x+vt,v), & {\rm if }\;  x+vt \leq l_j,
			\\
			0,& { \rm if\, not }.
		\end{cases}
	\end{align}
	for all $f\in X$, $(x,v)\in \Omega_j$ and $j\in\calM$. Moreover, 
	\begin{align*}
		(R(\mu,A)f)_j(x,v)=\int_{x}^{l_j}	e^{\int_x^{y}\tfrac{q_j(\sigma,v)-\mu}{v} d\sigma}\tfrac{1}{v}f_j(y,v)dy
	\end{align*}
	for $f\in X$, $(x,v)\in \Omega_j$, $\mu\in \rho(A)=\C$ and $j\in \calM$. Thus, by the additivity of the norm on the positive cone and Fubini’s theorem, we obtain
	\begin{align*}
		\Vert R(\mu,A)f\Vert_{X}&=\sum_{j=1}^{m}\int_{v_{\min}}^{v_{\max}}\int_{0}^{l_j} (R(\mu,A)f)_j(x,v) dxdv\\
		&= \sum_{j=1}^{m}\int_{v_{\min}}^{v_{\max}}\int_{0}^{l_j} \int_{x}^{l_j}	e^{\int_{x}^{y}\frac{q_j(\sigma,v)-\mu}{v} d\sigma}\frac{1}{v}f_j(y,v)dydxdv\\
		&= \sum_{j=1}^{m}\int_{v_{\min}}^{v_{\max}}\int_{0}^{l_j} \int_{x}^{l_j}	e^{-\int_{0}^{x}\frac{q_j(\sigma,v)-\mu}{v} d\sigma}e^{\int_{0}^{y}\frac{q_j(\sigma,v)-\mu}{v} d\sigma}\frac{1}{v}f_j(y,v)dydxdv\\
		&\ge \sum_{j=1}^{m}\int_{v_{\min}}^{v_{\max}}\int_{0}^{l_j} \left(\int_{0}^{y} e^{\tfrac{ (\mu-\tilde{q}) }{v}x }dx\right)	e^{-\mu\frac{ y}{v}}e^{\int_{0}^{y}\frac{q_j(\sigma,v)}{v} d\sigma}  \frac{1}{v}f_j(y,v)dydv\\
		&\ge  \sum_{j=1}^{m}\int_{v_{\min}}^{v_{\max}}\int_{0}^{l_j}     \left[\tfrac{1}{\mu-\tilde{q}}e^{\tfrac{ (\mu-\tilde{q})}{v}y} -\tfrac{1}{\mu-\tilde{q}} \right]e^{-\mu\frac{ y}{v}}e^{\kappa\frac{ y}{v} d\sigma}  f_j(y,v)dydv\\
		&\ge \frac{ 1}{\mu-\tilde{q}}\sum_{j=1}^{m}\int_{v_{\min}}^{v_{\max}}\int_{0}^{l_j}     \left[e^{\frac{ (\kappa-\tilde{q})}{v_{\min}}y}- e^{\frac{ (\kappa-\mu)}{v_{\min}}y} \right] f_j(y,v)dydv,
	\end{align*}
	for all $f\in X_+$ and $\mu>\tilde{q}:=\sup_{j\in \calM}  \Vert q_j\Vert_{\infty}$. On the other hand, for $ j \in \calM$, it follows from the Mean Value Theorem that there exists $\alpha_j\in (0,l_j)$ such that
	\begin{align*}
		\int_{0}^{1}   e^{-\mu\frac{ y}{v_{\max}} } f_j(y,v)dy=e^{-\mu\frac{ \alpha}{v_{\max}} } \int_{0}^{1}   f_j(y,v)dy, 
	\end{align*}
	for all $ 0\le f_j(\cdot,v)\in C(0,l_j)$ and a.e. $v\in [v_{\max},v_{\min}]$. Thus, by density,
	\begin{align}\label{inverse1}
		\Vert R(\mu,A)f\Vert_{X}\ge \frac{1}{\mu-\tilde{q}}\left[e^{\frac{ (\kappa-\tilde{q})}{v_{\min}}\alpha_m}- e^{\frac{ (\kappa-\mu)}{v_{\min}}\alpha_m} \right]\Vert f\Vert_{X}, \; \forall\,\mu>\bar{q},\;f\in X_+,
	\end{align}
	for $\alpha_m=\sup_{j\in \calM} \alpha_j$, where we have used the fact that $(\kappa-\tilde{q}) \ge (\kappa-\mu)$. On the other hand, the assumptions {\bf (A2)} and {\bf (A3)} imply that $G$ is surjective as $D(A_m)$ contains all constant functions $f$ satisfying 
	\begin{align*}
		f(l,v)=\mathcal{I}^{out}x,\qquad {\rm for \; some }\quad x\in \C^N.
	\end{align*}
	Moreover, a simple computation shows that the Dirichlet operator $D_\mu$ associated to $G$ is given by 
	\begin{align*}		
		(D_\mu g)_j(x,v)=e^{\int_{x}^{1}\tfrac{q_j(\sigma,v)-\mu}{v}d\sigma}\sum_{i\in \calN}\mathsf{w}_{ij}g_i(v),
	\end{align*}
	for $g\in \partial X$, $(x,v)\in \Omega_j$, $\mu\in \rho(A)$ and $j\in \calM$. Clearly, $D_\mu $ is positive and hence the assumptions (\emph{i}) and (\emph{ii}) are satisfied. 
	
	Let us now select 
	\begin{align*}
		C:=\Gamma_{\vert D(A)}, \qquad\qquad B:=(\mu I_X -A_{-1})D_\mu, \quad \mu \in \rho(A).
	\end{align*}
	It follows from Theorem \ref{S2.T1} and Lemma \ref{S2.P3} that $B$ and $C$ are L$^p$-admissible positive control and observation operators for $A$, since $A$ satisfies the inverse estimate \eqref{inverse1} for all $\mu>\tilde{q}$ and $\partial X$ is an L$^1$-space, respectively. In addition, using the uniqueness of Laplace transform, we show that the input-maps $\Phi_t:L^1_+(\R_+;\partial X)\to X_+ $ of $(A,B)$ are given by 
	\begin{align*}
		\left(	\Phi_t v\right)_j(x,v) =\sum_{i=1}^{N}	e^{\int_x^{l_j}\frac{q(\sigma,v)}{v}d\sigma}\mathsf{w}_{ji}g_i(\tfrac{tv-l_j+x}{v}) \1_{[\tfrac{l_j-x}{v},\infty)}(t),
	\end{align*}
	for $t\ge 0,\; (x,v)\in\Omega_j,\;j\in \calM$, and $v\in L^1_+(\R_+;\partial X)$. Then, for $\t>0$ and $v\in W^{1,1}([0,\tau],\partial X)$ with $v(0)=0$, it make sense to define the following operator
	\begin{align*}
		(\mathbb{F}v)(t)&=\Gamma\Phi_t^{A} v \cr
		&=\sum_{j=1}^{M}\sum_{i=1}^{N}	\textsf{i}_{kj}^{out}\int_{v_{\min}}^{v_{\max}}\ell_j(0,v,v')e^{\int_0^{l_j}\tfrac{q(\sigma,v')}{v'}d\sigma}\mathsf{w}_{ji}g_i(\tfrac{tv'-l_j}{v'}) \1_{[\frac{l_j}{v'},\infty)}(t)dv',
	\end{align*}
	for all $t\ge 0$, $i,k\in \calN$ and $j\in \calM$. With this explicit expression of the input-output map $\mathbb{F}$, it is not difficult to prove that $\F$ satisfies the estimate \eqref{S4.Fu} on the positive cone of $W^{1,1}_{0}([0,\t];\partial X)$ for every $\t \ge 0$. Thus, according to Proposition \ref{S4.T1}, $(A,B,C)$ is a positive $L^1$-well-posed triplet on $X,\partial X,\partial X $. Corollary \ref{S4.C2} further yields that $(A,B,C)$ is a positive $L^1$-well-posed strongly regular triplet. In addition, for $t<\inf_{j\in \calM} \tfrac{l_j}{v}$, we have 
	\begin{align*}
		I_{\partial X}-\mathbb{F}= I_{\partial X},
	\end{align*}
	and hence by virtue of Definition \ref{S4.D4} the identity $I_{\partial X}$ is a positive admissible feedback operator for $(A,B,C)$. It follows that $(A,B,C)$ is a positive $L^1$-well-posed strongly regular triplet on $X,\partial X,\partial X $ with feedthrough zero and the identity $I_{\partial X}$ as a positive admissible feedback operator. Therefore, according to Theorem \ref{S4.CB}, the transport networks systems $(\Sigma_{\mathsf{TN}})$ has a unique positive mild solution given by the variation of constant formula \eqref{mild-solution}. This completes the proof.	
\end{proof}

\section{Conclusion}
In this work, we investigated the well-posedness and positivity property of infinite-dimensional linear system with unbounded input and output operators. Indeed, we described the structural properties of $L^p$-admissible control and observation operators in the Banach lattice setting. This allowed us to characterize the internal and external positivity of input-output systems. Moreover, we proved that weak regularity and strong regularity of positive $L^p$-well-posed linear systems are equivalent. Then, we obtained that every positive $ L^1 $-well-posed system on Banach lattices is strongly regular. Furthermore, we provided sufficient conditions, in terms of an inverse estimate with respect to the Hille-Yosida Theorem, for zero-class admissibility for positive semigroups. We point out that, in practical situations, the state space is a $L^1$-space. Finally, we established tow generation results on positive perturbations of positive semigroups, namely the Desch-Schappacher perturbation and the Staffans-Weiss perturbation.

\end{document}